\documentclass[12pt]{article}\usepackage[]{graphicx}\usepackage[]{color}
\makeatletter
\def\maxwidth{ %
  \ifdim\Gin@nat@width>\linewidth
    \linewidth
  \else
    \Gin@nat@width
  \fi
}
\makeatother

\definecolor{fgcolor}{rgb}{0.345, 0.345, 0.345}

\usepackage{framed}
\makeatletter
 {\par\unskip\endMakeFramed%
 \at@end@of@kframe}
\makeatother

\definecolor{shadecolor}{rgb}{.97, .97, .97}
\definecolor{messagecolor}{rgb}{0, 0, 0}
\definecolor{warningcolor}{rgb}{1, 0, 1}
\definecolor{errorcolor}{rgb}{1, 0, 0}

\usepackage{alltt}
\usepackage{graphicx}
\usepackage{amssymb,amsmath,amsthm}
\usepackage{mathtools}
\usepackage{natbib,rotating}
\usepackage{bbm}
\usepackage{enumitem}
\usepackage{url}

\newtheorem{theorem}{Theorem}[section]
\newtheorem{definition}[theorem]{Definition}

\newtheorem*{remark}{Remark}

\newcommand{\Alpha}{\mathrm{A}}
\newcommand{\Beta}{\mathrm{B}}
\renewcommand{\P}{\mathbb{P}}
\newcommand{\R}{\mathbb{R}}

\newcommand{\Y}{\mathrm{Y}}
\newcommand{\U}{\mathrm{U}}
\renewcommand{\S}{\mathrm{S}}
\newcommand{\E}{\mathbb{E}}

\IfFileExists{upquote.sty}{\usepackage{upquote}}{}
\begin{document}

\title{Interpretable High-Dimensional Inference Via Score Projection with an Application in Neuroimaging}
 \author{ Simon N. Vandekar, Philip T. Reiss, and Russell T. Shinohara*}

  \maketitle
   \begin{center}
  \textbf{Authors' Footnote:}
\end{center}
Simon N. Vandekar is Doctoral Candidate (\url{simonv@mail.med.upenn.edu}),
Department of Biostatistics, Epidemiology, and Informatics, University of Pennsylvania, Philadelphia, PA 19104. Philip T. Reiss is Associate Professor (\url{reiss@stat.haifa.ac.il}), Department of Statistics, University of Haifa, Haifa, Israel.  Russell T. Shinohara is Assistant Professor of Biostatistics (\url{rshi@mail.med.upenn.edu}), Department of Biostatistics, Epidemiology, and Informatics, University of Pennsylvania, Philadelphia, PA 19104. SNV was supported by NIMH grant T32MH065218-11; PTR was supported by NIMH grant R01MH095836 and Israel Science Foundation grant 1777/16; RTS was supported by NINDS grant R01NS085211 and R21NS093349 as well as NIMH grant R01MH112847. The content is solely the responsibility of the authors and does not necessarily represent the official views of the National Institutes of Health. The authors thank Wei Pan and Haochang Shou for helpful discussions related to this work. Code to perform the analyses in this manuscript is provided at \url{https://bitbucket.org/simonvandekar/pst}. The data used in this manuscript are publicly available at \url{http://adni.loni.usc.edu/}.

~\\
*Corresponding author

%

\newpage
\begin{center}
\textbf{Abstract}
\end{center}
In the fields of neuroimaging and genetics, a key goal is testing the association of a single outcome with a very high-dimensional imaging or genetic variable.
Often, summary measures of the high-dimensional variable are created to sequentially test and localize the association with the outcome.
In some cases, the results for summary measures are significant, but subsequent tests used to localize differences are underpowered and do not identify regions associated with the outcome.
Here, we propose a generalization of Rao's score test based on projecting the score statistic onto a linear subspace of a high-dimensional parameter space.
In addition, we provide methods to localize signal in the high-dimensional space by projecting the scores to the subspace where the score test was performed.
This allows for inference in the high-dimensional space to be performed on the same degrees of freedom as the score test, effectively reducing the number of comparisons.
Simulation results demonstrate the test has competitive power relative to others commonly used.
 We illustrate the method by analyzing a subset of the Alzheimer's Disease Neuroimaging Initiative dataset.
Results suggest cortical thinning of the frontal and temporal lobes may be a useful biological marker of Alzheimer’s risk.


\newpage

\section{Introduction}
In scientific fields where high-dimensional data are prominent, significant interest lies in testing the association of a single continuous or categorical outcome with a large number of predictors. A common approach used in neuroimaging is to perform sequential tests to reduce the number of hypothesis tests. For example, it is common to first perform a test for the association of a phenotype with an imaging variable averaged across the entire brain. If the test rejects the null hypothesis of no association between brain and phenotype, then subsequent tests are conducted on regional averages of the data or on every voxel in the image using multiplicity correction to address the number of tests performed.
Often, location-specific results yield few or no significant findings due to reduced signal and the necessary adjustment for the large number of tests, even though the whole brain average data show a significant association.

In this paper, we propose a unified approach to test the association of an imaging or other high-dimensional predictor with an outcome and perform \emph{post hoc} inference to localize signal.
The framework is a modification of Rao's score test for models with a high-dimensional or infinite-dimensional parameter.
The theory developed assumes the parameter being tested is defined on a compact space such as the brain.
Though the approach is designed for hypothesis testing in neuroimaging, it is applicable to a wide range of scientific domains.


The standard framework assumes a model where $\mathrm{Y}_i$ are iid observations from density $f(y; \theta)$ and that the parameter $\theta = (\alpha, \beta) \in \Theta \subset \R^{m+p}$ where $\alpha \in \R^m$ is a nuisance parameter and $\beta \in \R^p $ is the parameter of interest.
We seek to test the hypothesis $H_0: \beta = \beta_0$.
Define the score function $\mathrm{U} = \mathrm{U}(\theta) = n^{-1} \sum_{i=1}^n\frac{\partial \log f(\mathrm{Y}_i \mid \theta) }{\partial \beta}(\theta)$.
Let $\theta_0 = (\alpha, \beta_0)$ be the null value of the parameter, where $\alpha$ is the true value of the parameter.
Let $\mathrm{S} = \U( \hat \alpha, \beta_0 )$ be the score function evaluated at the maximum likelihood estimate of $\alpha$ under the null hypothesis $H_0$.
Under the null and the conditions described in Appendix \ref{conditions}, the covariance of $\S$ can be obtained from the Fisher information evaluated at the null parameter value, 
\[
\Omega(\theta_0) = \mathbb{E} \left\{[(\partial/\partial \theta) \log f(\mathrm{Y}_1 \mid \theta)]^T [(\partial/\partial\theta) \log f(\mathrm{Y}_1 \mid \theta) ] \vert_{\theta_0} \right\}.
\]

The sum of scores (Sum) test originally discussed by \citet{rao_large_1948} has been used in genetics and neuroimaging \citep{pan_asymptotic_2009,kim_comparison_2014,madsen_groupwise_2009}.
The Sum test is based on the statistic
\begin{equation}\label{rao}
n\frac{(\mathrm{S}^T \zeta)^2}{\zeta^T \hat \Omega \zeta},
\end{equation}
where $\zeta \in \R^p$ is a given vector of weights. The denominator is an estimate of the variance of 
$\mathrm{S}^T \zeta$, so that the statistic is asymptotically $\chi^2_1$ \citep{rao_large_1948}.
The Sum test is locally most powerful \citep{rao_large_1948,cox_theoretical_1979}, however, the test has low power when there is a large number of variables that are not associated with the outcome \citep{pan_powerful_2014}.
This is due to the fact that the variance of the statistic in the numerator of \eqref{rao} increases by adding variables unassociated with the outcome without increasing the expected value of the numerator.

In the case of unknown weights, when $p<n$, \citet{rao_large_1948} proposed maximizing the Sum test statistic with respect to the weights,
\begin{equation}\label{raos}
\max_{\zeta\ne 0} n\frac{(\mathrm{S}^T \zeta)^2}{\zeta^T \hat \Omega \zeta} = n\mathrm{S}^T\hat\Omega^{-1} \mathrm{S}.
\end{equation}
This statistic is distributed as $\chi^2_p$ under the null. When $n >p$, \eqref{raos} is the usual score statistic, however this statistic cannot be used for high dimensional data due to the estimate $\hat \Omega$ being noninvertible when $p>n$.

For finite dimensional parameters, our proposed test can be thought of as a generalization of Rao's test in the case where the estimate of the information matrix is noninvertible.
When $p>n$, the test maximizes the statistic \eqref{rao} with respect to the vector $\zeta$ over a subspace, $\mathbb{L}$, of $\R^p$.
Maximization of the Sum test in the subspace $\mathbb{L}$ is equivalent to projecting the scores for the original model to a lower dimensional space where the information matrix is invertible.
For this reason, we call the test a projected score test (PST).
The procedure does not assume sparsity, but attempts to conserve power by reducing the dimension of the data and performing inference in the lower dimensional space.

In many cases, if a score test rejects $H_0$, then it is of primary interest to perform \emph{post hoc} inference to identify nonzero parameters.
In neuroimaging, this amounts to a high-dimensional testing problem where the association is tested at each location in the image.
The standard approach is to perform a hypothesis test at each parameter location and use a multiplicity correction procedure.
Such methods this in neuroimaging that control the family-wise error rate (FWER) have relied on Gaussian random field theory \citep{friston_assessing_1994}, but have recently been shown to have type 1 error rates far from the nominal level in real data due to assumption violations \citep{eklund_cluster_2016,silver_false_2011}.
Recently, considerable research activity has focused on leveraging the dependence of the tests to control the false discovery rate (FDR) in high-dimensional settings \citep{efron_correlation_2007}.
\citet{sun_false_2015} develop a procedure to control the FDR for spatial data as well as an approach for controlling the expected proportion of false clusters.
\citet{fan_estimating_2012} discuss estimation of the false discovery proportion (FDP) under dependence for normally distributed test statistics based on a factor approximation.
In contrast, the PST \emph{post hoc} inference procedure is performed by projecting the scores onto $\mathbb{L}$, and controlling the FWER of the projected scores.



Several recent studies have considered hypothesis tests for functional data, which is conceptually similar to our approach for an infinite-dimensional parameter.
\citet{reiss_functional_2010} propose inverting simultaneous confidence bands for the parameter of a functional predictor to test which locations of the image are associated with the outcome.
\citet{smith_spatial_2007} use a binary Markov random field model to compute the joint probability that the marginal parameter estimates are equal to zero.
Our \emph{post hoc} inference is most similar to \citet{smith_spatial_2007} as the interpretation of the contribution of the scores retains a marginal interpretation.

Here, we derive the asymptotic null distribution for the PST statistic under some standard regularity conditions.
For data that are measured on a compact space, such as brain images, we discuss sufficient theoretical assumptions for characterizing test behaviors as both $n$ and $p$ approach infinity.
For a normal linear model, we show how the finite sample distribution of the statistic can be calculated exactly for fixed $n$ and $p$.

To demonstrate how the test can be used in neuroimaging, we investigate the association of cortical thickness with mild cognitive impairment (MCI) in the Alzheimer's Disease Neuroimaging Initiative (ADNI) study, a data set where $p=$18,715 and $n=628$.
The outer surface of the brain (cortex) represents a highly folded sheet in 3-dimensional space.
The thickness of the cortex is known to be affected in individuals with psychopathology and neurological illness.
MCI is a subtle pre-Alzheimer's disease decline in cognitive functioning.
There is significant clinical interest in finding biological markers of MCI in order to identify those at risk for developing Alzheimer's disease, as prevention strategies and therapeutics for early disease are increasingly common. 
In this data set, we seek to localize regions of the brain where cortical thinning provides additional information with regard to the diagnosis of MCI beyond what can be ascertained by neurocognitive scales alone.


For the remainder of the manuscript, we denote matrices by uppercase italic letters ($X$), vectors by lowercase ($x$), and random vectors by uppercase roman letters ($\mathrm{X}$).
Hilbert spaces are denoted with black-board letters ($\mathbb{X}$) and Greek letters denote model parameters.
For the singular value decomposition (SVD) of any matrix we will assume that the smallest dimensions of the matrices obtained are equal to the rank of the matrix $X$ unless otherwise noted.
$\xrightarrow{L}$ denotes convergence in law and $\xrightarrow{P}$ denotes convergence in probability.



 \section{The Projected Score Test}
In Section \ref{sec:MaxSumExponential} and Appendix \ref{sec:framework} we define the PST statistic, give its asymptotic distribution, and lay out the theoretical framework.
In Section \ref{sec:asympp}, we detail conditions sufficient for studying asymptotics in $p$.
We discuss maximization of the sum test for normal linear models in Section \ref{sec:MaxSumNormal}.

\subsection{PST for finite-dimensional parameters}
\label{sec:MaxSumExponential}

We assume the observed data are finite-dimensional representations that are generated from an underlying stochastic process.
In Appendix \ref{sec:framework} we describe how to define the finite-dimensional likelihood from the infinite-dimensional likelihood.
Here, we informally define the finite-dimensional likelihood, the interested reader can refer to Appendix \ref{sec:framework} for further details.

Let $\mathbb{V}$ be a nonempty compact subset of $\R^3$ and $\mathcal{B}(\mathbb{V})$ be the space of square integrable functions from $\mathbb{V}$ to $\R$.
$\mathbb{V}$ represents the space on which data can be observed; in neuroimaging this space is the volume of the brain.
The underlying stochastic processes are assumed to take values in $\mathcal{B}(\mathbb{V})$, but the observed finite-dimensional data are $p$-dimensional discretizations of the stochastic processes defined on $\mathbb{V}$.
Thus, the observed data can be described as iid observations $\Y_i = (\Y_i^{(1)}, \Y_i^{(2)})$, for $i=1,\ldots, n$, with $\Y_i^{(1)}$ taking values in $\R^k$, and $\Y_{ip}^{(2)} \in \R^p$.
Observations of $\Y_i$ are a vector of $k$ variables that are nonimaging covariates and the outcome variable, together with the observed finite-dimensional neuroimaging data.
We denote the collection of data by $\Y = (\Y_1, \ldots, \Y_n)$.
We define a parameter space $\Theta = \R^m \times \R^p$ that includes a finite-dimensional nuisance parameter $\alpha \in \R^{m}$ and the $p$-dimensional discretized parameter of interest $\beta_p \in \R^p$.
Together these parameters describe the joint distribution of the imaging and nonimaging data.

Denote the finite dimensional likelihood by $\ell(\theta_p; \Y)$, where $\theta_p = (\alpha, \beta_p)$ and $ \Y$ are the discretized parameters and data, respectively. 
Define the score function $\U_{np} = \frac{\partial \ell}{\partial \beta_p}\{\theta_p; \Y\}$ and let 
\begin{equation}
\label{eq:frechetfinitederiv}
\mathrm{S}_{np} = \mathrm{U}_{np}(\hat\alpha, \beta_{p0})  \in \R^p
\end{equation}
denote the score function evaluated at the maximum likelihood estimate (MLE) under the null hypothesis $H_0: \beta_p = \beta_{p0}$.

Let the Fisher information for the full model be
\begin{equation}
\label{eq:info}
\Omega_F(\theta_{p0}) = 
\E_{\theta_{p0}} \left(\{\frac{\partial}{\partial \theta_p}\log f(\Y_1 \mid \theta_p)\}^T \{\frac{\partial}{\partial \theta_p}\log f(\Y_1 \mid \theta_p)\} \Big\vert_{\theta_{p0}} \right) = 
\begin{bmatrix}
\Omega_{\alpha} & \Omega_{\alpha\beta} \\
\Omega_{\beta\alpha} & \Omega_{\beta}
\end{bmatrix},
\end{equation}
where $\theta_{p0} = (\alpha, \beta_{p0})$.
Then the asymptotic variance for $\sqrt{n}\mathrm{S}_{np}$  under $H_0$ is
\begin{equation}\label{eq:effinfo}
\Omega(\theta_{p0}) = \Omega_\beta - \Omega_{\beta\alpha}\Omega^{-1}_{\alpha} \Omega_{\alpha\beta}.
\end{equation}
With the finite parameter scores defined, we can define the PST.


\begin{definition}\label{def:MaxSum}
Let $P_\mathbb{L}$ be the orthogonal projection matrix onto a linear space $\mathbb{L} \subset \R^p$ with $r = \text{dim}(\mathbb{L}) < n-m$.
Let $\mathrm{S}_{np}$ be as defined in \eqref{eq:frechetfinitederiv}
and $\hat\Omega$ be the plug-in estimator of the covariance \eqref{eq:effinfo} obtained from 
\begin{equation}\label{infoestimate}
\hat \Omega_F = n^{-1}\sum_{i=1}^n \left(\frac{\partial\log f(\Y_i ; \theta_p)}{\partial\theta_p^T} \right)\left(\frac{\partial\log f(\Y_i \mid \theta_p)}{\partial\theta_p} \right) \Bigr\vert_{\hat\theta_{p0}},
\end{equation}
where $\hat \theta_{p0} = (\hat \alpha, \beta_{p0})$ denotes the maximum likelihood estimate of the parameter vector under the null hypothesis $H_0 : \beta = \beta_{p0}$.
Then the PST statistic with respect to $\mathbb{L}$ is defined as
\[
\mathrm{R}^{\mathbb{L}}
= \max_{\zeta \in \mathbb{L}\setminus\{0\}} n\frac{(\zeta^T \mathrm{S}_{np})^2}{\zeta^T \hat\Omega(\theta_{p0}) \zeta }
 = \max_{\zeta \in \R^p\setminus\{0\}} n\frac{(\zeta^T P_{\mathbb{L}} \mathrm{S}_{np})^2}{\zeta^T P_{\mathbb{L}} \hat \Omega(\theta_{p0}) P_{\mathbb{L}} \zeta }.
\]
\end{definition}
The following theorem states that the asymptotic distribution (with respect to $n$) of the PST statistic can be found for any finite dimensional likelihood based on independent observations provided the same regularity conditions required for the convergence of the scores to a multivariate normal random variable.

\begin{theorem}\label{thm:MaxSumGeneral}
Assume all objects are as described in Definition \ref{def:MaxSum}.
Let $P_{\mathbb{L}} = Q Q^T$ where the columns of the $r\times p$ matrix $Q$ are any orthonormal basis for $\mathbb{L}$.
Define 
\[
V = V(\theta_{p0} ) = Q^T \Omega(\theta_{p0}) Q,
\]
and assume the estimate $\hat{\mathrm{V}} = Q^T \hat \Omega(\theta_{p0}) Q$ is invertible, and that the conditions given in Appendix \ref{conditions} are satisfied.

Then, under the null, the rotated scores, denoted $\S^Q_{np}$, are
\begin{equation}\label{eq:rotatedscores}
n^{1/2}\S^Q_{np} = n^{1/2} Q^T \mathrm{S}_{np} \xrightarrow{L} S^Q_p \sim N_r( 0 , V),
\end{equation}
the PST statistic is
\begin{equation}\label{thm1maxsum}
\mathrm{R}^{\mathbb{L}} = n (\mathrm{S}_{np}^Q)^T \hat{\mathrm{V}}^{-1} \mathrm{S}^Q_{np},
\end{equation}
and $\mathrm{R}^{\mathbb{L}} \xrightarrow{L} \chi^2_{r}$ as $n\to \infty$.
\end{theorem}

Theorem \ref{thm:MaxSumGeneral} requires that $\hat{\mathrm{V}}$ is nonsingular, however, in practice it is possible to ensure that $Q$ is in the column space of $\hat \Omega(\theta_{p0})$, so that $\hat{\mathrm{V}}^{-1}$ exists.
The proof of Theorem \ref{thm:MaxSumGeneral} is given in Appendix \ref{prooftheorem1}.
We also demonstrate there that the result of Theorem \ref{thm:MaxSumGeneral} does not depend on the choice of $Q$.
We show how $\mathbb{L}$ can be chosen for GLMs in Sections \ref{sec:specialcases} and \ref{sec:adaptPCA} and for imaging data in the analysis of the ADNI dataset in Section \ref{sec:DataAnalysis}.

\subsection{The PST as $p\to \infty$}
\label{sec:asympp}

We will show that as $p\to\infty$ the PST statistic converges to an integral over a stochastic process.
The rate that $p$ approaches infinity does not depend on the sample size.
Here, we assume the data can take values on the Hilbert space $\mathbb{Y} =\R^k\times \mathcal{B}(\mathbb{V}) $, where $\mathbb{V}$ is a nonempty compact subset of $\R^3$ and $\mathcal{B}(\mathbb{V})$ is the space of square integrable functions from $\mathbb{V}$ to $\R$.
Let $\Y_i = (\Y_i^{(1)}, \Y_i^{(2)})$, for $i=1,\ldots, n$, be iid with $\Y_i^{(1)}$ taking values in $\R^k$, and $\Y_i^{(2)}$ a stochastic process taking values in $\mathcal{B}(\mathbb{V})$.
Realizations of $\Y_i$ are a vector of $k$ variables that are nonimaging covariates and the outcome variable, together with a function on $\mathbb{V}$.
We assume the parameter $\beta \in  \mathcal{B}(\mathbb{V})$.
The infinite-dimensional score function is defined in Appendix \ref{sec:framework} as the Fr\'echet derivative of the likelihood with respect to the parameter $\beta$, $\U_n = \U_n(v) = \frac{\partial \ell}{\partial \beta}\{(\alpha, \beta(v)); \Y(v)\}$.
And the score is defined as the function
\begin{equation}
\label{eq:infinitederiv}
\mathrm{S}_{n} = \mathrm{U}_{n}\{\cdot ; (\hat\alpha, \beta_0) \} \in \mathcal{B}(\mathbb{V}).
\end{equation}
Throughout we assume that the infinite-dimensional scores converge in law, i.e.
\begin{equation}\label{eq:scorefuncconvergence}
n^{1/2}\mathrm{S}_n \to_L \mathrm{S},
\end{equation}
where $\mathrm{S}$ is a mean zero Gaussian process.
Theorem \ref{thm:stochasticprocessCLT} in Appendix \ref{sec:extras_asympp} gives conditions under which this convergence holds \citep{van_der_vaart_asymptotic_2000}.
The following definition of the PST statistic for infinite-dimensional parameters is motivated by formula \eqref{thm1maxsum}.

\begin{definition}
\label{def:MaxSumInf}

Let $(q_1(v), \ldots, q_r(v))$ be an orthonormal basis for the linear subspace $\mathbb{L} \subset \mathcal{B}(\mathbb{V})$ with respect to the $L^2(\nu)$ inner product where $\nu$ is the Lebesgue measure, and $r = \text{dim}(\mathbb{L}) < n-m$.
Assume $q_j$ are such that
\begin{equation}\label{eq:discontinuities}
\nu( \{ v : \text{$q_j$ is discontinuous at v} \} )=0
\end{equation}
for all $j=1,\ldots,r$.

Define the column vector $\S^Q_n \in \R^r$, with $j$th element 
\begin{equation}\label{eq:intinnerproduct}
(\S^Q_n)_j = \int_{\mathbb{V}} q_j(v) \S_n(v)dv,
\end{equation}
and let $\hat{ \mathrm{V}}_n$ be the $r \times r$ covariance matrix with $(j,k)$th element
\begin{align*}
\hat{ \mathrm{V}}_n^{j,k} = & n^{-1}\sum_{i=1}^n
\left(\int_\mathbb{V} q_j(v)\left[\frac{\partial}{\partial \beta}\log f\{\Y_i(v) \mid \hat\alpha, \beta_0(v)\}\right]dv \right) \\
& \times \left(\int_\mathbb{V} q_k(v)\left[\frac{\partial}{\partial \beta}\log f\{\Y_i(v) \mid \hat\alpha, \beta_0(v)\}\right]dv \right),
\end{align*}
where $\frac{\partial}{\partial \beta}\log f\{\Y_i(v) \mid \hat\alpha, \beta_0(v)\}$ denotes the Fr\'echet derivative evaluated at $\beta_0$.
Let $V_n = \E \hat{ \mathrm{V}}_n $.
Assume that $\hat{ \mathrm{V}}_n$ is invertible. 
Then the PST statistic with respect to $\mathbb{L}$ is defined as
\[
\mathrm{R}^{\mathbb{L}}
= n(\S^Q_n)^T \hat{ \mathrm{V}}_n^{-1} \S^Q_n.
 \]

\end{definition}

While we have given a definition of the PST statistic in infinite dimensions, in practice this statistic is not estimable because it depends on functions which are only observed on finite a grid. The following theorem states that as the resolution of the grid is increased then the finite parameter PST statistic converges to the PST statistic for the infinite-dimensional parameter.
Moreover, as the sample size increases the statistic converges to a statistic based on the Gaussian process $\S$ \eqref{eq:scorefuncconvergence}.
The rate that $p$ increases does not depend on $n$.

\begin{theorem}\label{thm:maxsumconvergence}
Let $\S_{np}$ be as defined in \eqref{eq:frechetfinitederiv}.
For objects as defined in Definition \ref{def:MaxSumInf}, let $q_{jp} = (q_j(v_{1p})\nu(\mathbb{V}_{1p}), \ldots, q_j(v_{pp})\nu(\mathbb{V}_{pp}) )^T$, where $v_{jp}$ and $\mathbb{V}_{jp}$ are as defined in Appendix \ref{sec:framework}. 
Let $Q_p$ be the $p \times r$ matrix with $j$th column $q_{jp}$.
Denote $\S^Q_{np} = Q^T_p \S_n$.
Define the $j$th element of the vector $\S^Q \in \R^r$ as
\begin{equation*}
(\S^Q)_j = \int_{\mathbb{V}} q_j(v) \S(v)dv,
\end{equation*}
Assume the conditions for Theorems \ref{thm:MaxSumGeneral} and \ref{thm:stochasticprocessCLT}, and that $\S_n$ and $\frac{\partial}{\partial \beta}\log f\{\Y_i(v) \mid \hat\alpha, \beta_0(v)\}$ have continuous sample paths with respect to $v$ (i.e. for almost every $\Y \in \mathbb{Y}^{n}$ and $\theta \in \Theta$, $\S_n(\cdot;\Y, \Theta)$ is continuous).
Let $V = \lim_{n\to \infty} \E V_n$.
For $p_1 > p_2$, let $\mathcal{V}_{p_1}$ be a refinement of $\mathcal{V}_{p_2}$ such that 
\begin{equation}\label{eq:intconvergence}
\lim_{p\to\infty} \sup_{k} \nu(\mathbb{V}_{kp}) = 0.
\end{equation}
Then as $n,p \to \infty$,
\begin{equation}\label{eq:maxsumconvergence}
n (\S^Q_{np})^T \hat{\mathrm{V}}_{np}^{-1} \S^Q_{np}
\xrightarrow{P} (\S^Q)^T V^{-1} \S^Q.
\end{equation}
\end{theorem}
The proof is given in Appendix \ref{sec:proofasympp}.

\subsection{The PST in Normal Linear Models}
 \label{sec:MaxSumNormal}

The finite-sample distribution for the PST statistic can be found exactly for a normal linear model.
Define $X = [x_1, \ldots, x_n]^T$ to be an $n\times m$ full rank matrix of covariates for each observation, $\tilde G = [g_1, \ldots, g_n]^T$ to be an $n \times p$ full rank matrix of predictor variables of interest with $p>n$, and $\tilde \Y = [ \tilde \Y_1, \ldots, \tilde \Y_n ]^T$ to be $n \times 1$ normal random vector with independent elements conditional on $X$ and $G$.
 The Sum test with normal error is based on the model
 \[
 \tilde \Y_i = \alpha^T x_i + \beta^T g_i + \mathrm{E}_i,
 \]
 where all variables are as previously defined and $\mathrm{E}_i \sim N(0, \sigma^2)$ are independent.
 If we let $A A^T =  (I-H)$ be the SVD of the projection $(I-H)$, where $H = X(X^TX)^{-1}X^T$, and define $G = A^T \tilde G$ and $Y = A^T \tilde \Y $, then under the null $Y  \sim N_{(n-m)}(0, \sigma^2 I)$.
 
 The Sum test statistic for $H_0 : \beta = 0$ is \citep{rao_large_1948,lin_general_2011}
 \begin{equation}\label{normalsum}
 \mathrm{R}_{Sum} = n\frac{(\mathrm{S}_{np}^T \zeta)^2}{\zeta^T \hat \Omega \zeta} = \frac{(\Y^T G \zeta)^2}{ \hat\sigma^2 \zeta^T G^T G \zeta },
 \end{equation}
 where $\zeta \in \R^p$ is a known vector of weights, $\mathrm{S}_{np} = n^{-1} G^T \Y$  is the score vector evaluated at the MLEs under the null ($\mathbb{E}(\tilde \Y_i) =\alpha^Tx_i $), and $\hat \Omega = n^{-1}\hat\sigma^2 G^T G$ is an estimate of the variance of $\S_{np}$ which corresponds to using the estimator \eqref{infoestimate}.
 
%

 
 The PST statistic for this model with respect to a linear subspace $\mathbb{L}$ of $\R^p$ by Definition \ref{def:MaxSum} is
 \begin{align*}
 \mathrm{R}^{\mathbb{L}} 
 &= \max_{\zeta \in \mathbb{L}} \frac{(\Y^T G \zeta)^2}{ \hat\sigma^2 \zeta^T G^TG \zeta }.
 \end{align*}
The following theorem gives a closed-form expression for $\mathrm{R}^{\mathbb{L}}$ and its null distribution.

\begin{theorem}\label{thm:MaxSum}
Define $W$ be the projection matrix onto the column space of $G P_{\mathbb{L}}$.
Then
\begin{equation} \label{normalmaxsum}
\mathrm{R}^{\mathbb{L}} = (n-m)\frac{\Y^T W \Y}{\Y^T \Y},
\end{equation}
and under the null, 
\begin{equation}
 \label{eq:normalmaxsumdist}
\mathrm{R}^{\mathbb{L}} =_L \frac{r(n-m)}{r + (n-m -r) \mathrm{F}_{(n-m-r),r} },
\end{equation}
where $\mathrm{F}_{(n-m-r),r}$ is F-distributed with $(n-m-r)$ and $r$ degrees of freedom.
\end{theorem}
The proof can be found in Appendix \ref{prooftheorem2}. The form of equation \eqref{normalmaxsum} shows that for a normal linear model, the test statistic is a ratio of quadratic forms.
Due to the rotation invariance of $\Y$ under the null, the finite-sample distribution of $\mathrm{R}^{\mathbb{L}}$ depends only on the sample size and the dimension of the basis, but not on the particular choice of $\mathbb{L}$.



\section{Specifying the linear subspace $\mathbb{L}$}

\subsection{Specifying $\mathbb{L}$ in generalized linear models}
\label{sec:specialcases}

Here, we discuss choices for the selection of $\mathbb{L}$ in the context of GLMs with the canonical link function.
We restrict attention to finite dimensional parameters and forgo the subscripts on the finite sample score vector $\S$.
Define $X = [x_1, \ldots, x_n]^T$, and $G = [g_1, \ldots, g_n]^T$ where objects are as defined in Section \ref{sec:MaxSumNormal}.
Assume the outcome $\Y = [ \Y_1, \ldots, \Y_n ]^T$ is from an exponential family where the expectation can be written
\[
h(\E \Y_i) = \alpha^T x_i + \beta^T g_i,
\]
where $h$ is the canonical link function.
For the GLM with canonical link, the scores are \citep{mccullagh_generalized_1989} 
\begin{equation*}
\mathrm{S} = n^{-1} G^T (\Y - \hat \Y),
\end{equation*}
where 
\begin{equation}\label{eq:Yhat}
\hat {\Y} = [ \hat{\Y}_1, \ldots \hat{\Y}_n ]^T
\end{equation}
and $\hat {\Y}_i = h^{-1}( x^T_i \hat \alpha)$ is the $i$th fitted value under the null. 
Let $\Gamma$ be the $n\times n$ diagonal matrix with $i$th diagonal element $\Gamma_{ii} = (\Y_i - \hat{\Y}_i)^2$.
Then the estimate of the covariance \eqref{eq:effinfo} obtained using \eqref{infoestimate} is
\begin{equation*}
\hat \Omega = n^{-1}\{ G^T \Gamma G - G^T \Gamma X ( X^T \Gamma X)^{-1} X^T\Gamma G \}
\end{equation*}
The score statistic is obtained from the scores and the estimated information as in expression \eqref{raos}.

In this setup, the basis for $\mathbb{L}$ can be constructed from the principal component analysis (PCA) of $G$.
We write the PCA of $G$ in terms of the SVD $G = T_* D Q^T$, where the principal scores are $T = T_*D = GQ$.

With this basis, the PST is equivalent to performing Rao's score test in a principal components regression model.
To see this, first note that principal component regression is defined by
\[
h(\E\Y ) = X \alpha + T \beta_T.
\]
The scores for $\beta_T$ are
\[
\mathrm{S}_{T} = n^{-1} Q^T G^T (\Y - \hat \Y) = Q^T \mathrm{S},
\]
which are the same as the rotated scores in \eqref{eq:rotatedscores}.
The information estimate is also equivalent.
Thus the score test statistic, $n \mathrm{S}_T^T\hat \Omega_T^{-1} \mathrm{S}_T$, in principal component regression is equivalent to the PST statistic \eqref{thm1maxsum}.

Another useful basis for $\mathbb{L}$ may be constructed from vectors that are indicators of variables that are expected to have a similar relationship with the outcome.
The anatomical basis used in Section \ref{sec:DataAnalysis} is an example.
To define the basis vectors $q_j$, $j=1,\ldots, r$, we let $\mathcal{Q}_j \subset \{1, \ldots, p\}$ such that $\mathcal{Q}_j \cap \mathcal{Q}_{j'} = \varnothing$ for $j\ne j'$, and then set the $k$th element of the $j$th basis vector to be $q_{kj} = \mathbbm{1}( k \in \mathcal{Q}_j)$.
These define orthogonal basis vectors since the sets $\mathcal{Q}_j$ are disjoint.
This basis is equivalent to averaging $r$ subsets of the $p$ predictor variables and performing a hypothesis test of the regression onto the $r$ averaged variables.

The choice of the basis is a critical decision as it affects the power and interpretation of the \emph{post hoc} inference.
To clarify, under the alternative the scores have nonzero mean 
\begin{equation}\label{eq:scoremean}
\E \S = \mu \in \R^{p}.
\end{equation}
If the projection is orthogonal to $\mu$ then the test will have power equal to the type 1 error rate.
The PCA basis assumes that $\mu$ has a spatial pattern similar to the covariance structure of the predictor variables.
The anatomical basis assumes that all locations within a region have the same parameter value.
We discuss the effect of the basis on the interpretation of the \emph{post hoc} inference in Section \ref{sec:posthoc}.


\subsection{Choosing a dimension for the PCA basis}
\label{sec:adaptPCA}

In order to choose a dimension for the PCA basis, we propose an adaptive procedure that sequentially tests bases of increasing dimension while controlling the type 1 error rate.
To do this we first condition on the parameter estimate $\hat\alpha$ for the reduced model and perform the SVD $(\Gamma - \Gamma X (X^T \Gamma X)^{-1} X^T \Gamma)^{1/2}G = T D Q^T$.
We use subsets of columns $Q$ as the basis $\mathbb{L}$.
For two columns $q_j$ and $q_k$ of $Q$
\[
\text{Cov}(q^T_j \S, q_{k}^T \S) = q^T_j G^T (\Gamma - \Gamma X (X^T \Gamma X)^{-1} X^T \Gamma) G q_{k}
= q^T_j Q D^2 Q^T q_k= 0.
\]
Thus each projected score $n^{1/2} q^T_j \S$ is asymptotically independent because $n^{1/2}Q^T S$ is asymptotically normal and can be tested by a separate chi-squared test at level $\alpha^*$.
If this is done sequentially for $r=1,\ldots,n-m$, then, due to their asymptotic independence, the probability of a type 1 error under the global null is
\begin{align*}
\sum_{r=1}^{n-m}\P\left((q^T_j \S)^2 >\chi^2_1(\alpha^*) \text{ for all $j\le r$} \right)
& \approx \sum_{r=1}^{n-m} (\alpha^*)^{r} \\
& \le \sum_{r=1}^{\infty} (\alpha^*)^{r} = \frac{1}{1-\alpha^*} - 1,
\end{align*}
where $\chi^2_1(\alpha^*)$ denotes the $1-\alpha^*$ quantile of a chi-squared distribution. The approximate equality is due to the asymptotic approximation and the final inequality is the geometric series solution.
In order to control the type 1 error at level $\alpha$, we choose $\alpha^*= 1-(1+\alpha)^{-1}$.
Then we sequentially test $r=1,\ldots,(n-m)$ until we fail to reject a test at level $\alpha^*$.
Note that the power depends critically on the first test in the sequence; subsequent tests serve only to increase the dimension of the basis.
If the first component is orthogonal $\mu$ in \eqref{eq:scoremean}, the probability of reaching other components that do is less than $\alpha^*$.

A potentially more robust procedure is to test chunks of PCs by varying $r = \{r_1 = 0, r_2, r_3, \ldots, r_k = n-m\}$ and for the $j$th test perform a chi-squared test of all PCs $(r_{j}+1), \ldots,r_{j+1}$ on $r_{j+1} - r_{j}$ degrees of freedom.
So long as the tests are independent, which they are under the global null, the rejection threshold $\alpha^*$ will control the type 1 error rate at level $\alpha$.
This adaptive PCA (aPCA) procedure is implemented below by testing the first 5 components together and sequentially testing components $6,\ldots,n-m$ one-at-a-time.
We demonstrate the procedure in the ADNI data analysis below and type 1 error rates are assessed in Section \ref{sec:Simulations}.

\section{\emph{Post hoc} Inference for Localizing Signal}
\label{sec:posthoc}

After performing the test of association using the PST statistic, it is of primary interest to investigate the contribution of the scores to the statistic in order to identify which locations in the image are associated with the outcome and the direction of the effect.
This can be done by projecting the scores onto $\mathbb{L}$ and performing inference that controls the FWER for the projected scores.
Because the projected scores are distributed in a subspace of $\R^p$, inference is much less conservative compared to performing inference on the original score vector.




Our aim is to construct a rejection region for each element of the projected score vector $(P_{\mathbb{L}} \mathrm{S})_{j}$, for $j=1, \ldots p$. Under the null,
\[
P_{\mathbb{L}} \mathrm{S} \sim N(0, P_{\mathbb{L}}\Omega P_{\mathbb{L}}).
\]
The diagonal elements of $P_{\mathbb{L}}\Omega P_{\mathbb{L}}^T$ are not equal, so defining a single rejection threshold for all elements favors rejection for elements with larger variances.
To resolve this issue we scale by the inverse of the standard deviation of the projected scores.
Let $\Delta$ be the diagonal matrix with $j$th diagonal element $\Delta_{jj} = 1/\sqrt{(P_{\mathbb{L}} \Omega P_{\mathbb{L}})_{jj} }$.
Then the rejection threshold that controls the FWER for the standardized projected scores is defined by $c$ that satisfies
\begin{align}\label{eq:RR}
1-\P(\lvert (\Delta P_{\mathbb{L}}\mathrm{S})_j \rvert > c \text{ for any $j$}) 
& = \mathbb{P}( \max_j\lvert (\Delta P_{\mathbb{L}} \mathrm{S})_j \rvert < c) = 1 - \alpha.
\end{align}

Thus, the distribution of the infinity norm of $\Delta P_{\mathbb{L}}\S$ can be used to compute a rejection threshold for the standardized projected scores that controls the FWER.
The rejection threshold that controls the FWER is $c\in \R^+$ such that
\begin{align*}
 \mathbb{P}( \lvert \Delta P_{\mathbb{L}} \mathrm{S} \rvert_\infty < c) & = 1 - \alpha.
\end{align*}
This set defines the region where the probability any element of the standardized projected score vector is greater than $c$ is equal to $\alpha$ under the global null $H_0: \beta = \beta_0$.
We reject the null hypothesis at location $j$ if the observed projected score $\lvert (\Delta P_{\mathbb{L}} s)_j\rvert > c$.
This threshold corresponds to a single-step ``maxT" joint multiple testing procedure \citep{dudoit_multiple_2008} and satisfies the assumption of subset pivotality, so it controls the FWER at the nominal level in the case that some projected scores have nonzero mean \citep{westfall_resampling-based_1993} (see Appendix \ref{sec:subsetpivotality}).

By \eqref{eq:rotatedscores} we have
\[
\Delta P_{\mathbb{L}} \mathrm{S} \xrightarrow{L} \Delta Q V^{1/2} \mathrm{Z},
\]
where $\mathrm{Z} \sim N_r(0, I)$. Thus we can approximate the region in \eqref{eq:RR} by finding $c$ so that
\[
\int_{\lvert \Delta Q V^{1/2} z \rvert_\infty \le c} \phi_r(z) dz = 1 - \alpha,
\]
where $\phi_r$ denotes the PDF of $Z$.
In practice we approximate this interval by plugging in estimates for $\Delta$ and $V^{1/2}$.

This integral is difficult to calculate due to the large dimensions of $Q$, but can be approximated quickly and easily using Monte Carlo simulations.
$B$ simulations are used to estimate the CDF of the infinity norm, $\hat F_B(\cdot)$, which we use to obtain p-values for each observed standardized projected score, $(\Delta Ps)_j$, by evaluating 
\begin{equation}\label{eq:pvalues}
p_{j} = 1- \hat F_B\left((\Delta P_{\mathbb{L}}s)_j\right),
\end{equation}
or a rejection threshold can be obtained by using
\begin{equation}\label{eq:qvalues}
c = \hat F_B^{-1}(1-\alpha).
\end{equation}
The p-value for a given element of the standardized projected score vector is the probability of observing a projected score as large as $(\Delta P_{\mathbb{L}} s)_j$ under the global null $H_0: \beta = \beta_0$.
The standard deviation of the Monte Carlo estimate \eqref{eq:pvalues} decreases at a $\sqrt{B}$ rate and depends only on the volume of the space being integrated, so the procedure will perform well for computing adjusted p-values with a small error \citep{press_numerical_2007}.
For example, because the volume of the space being integrated is 1, with 10,000 simulations the standard deviation is on the order of $B^{-1/2} = 0.01$.

Rejection of the null hypothesis $H_0: \beta = \beta_0$ is not strictly necessary to proceed with the {\it post hoc} inference procedure; the {\it post hoc} procedure can be used separately from the PST.
In addition, it is important to note that the \emph{post hoc} inference can have improved power by interpreting the projected scores.
When the alternative hypothesis is true, the rejection regions for the projected scores do not necessarily control the type 1 error for the unprojected scores.
This is demonstrated in the imaging simulations in Section \ref{sec:Simulations}.

As mentioned above, the basis affects the interpretation of the inference on the projected scores.
For the PCA basis the interpretation is as follows: over repeated experiments if the data are projected onto $\mathbb{L}$, then the probability of falsely rejecting one or more scores $j$ with $(\Delta P_{\mathbb{L}} \mu)_j = 0$ is at most $\alpha$, where $\mu$ is as defined in \eqref{eq:scoremean}.

\section{ADNI Neuroimaging Data Analysis}
\label{sec:DataAnalysis}

We obtained data  from  the  Alzheimer's  Disease Neuroimaging  Initiative  (ADNI)  database (adni.loni.usc.edu).
The ADNI  was  launched  in 2003 as  a  public-private  partnership,  led  by  Principal  Investigator Michael  W.  Weiner, MD.
The ADNI is a longitudinal observational study designed to investigate the early biomarkers of Alzheimer's disease; detailed MRI methods are given by \citet{jack_alzheimers_2008}. 
Mild cognitive impairment (MCI) represents a subtle pre-Alzheimer's Disease decline in cognitive performance.
The goal of our analysis is to identify whether a subset of the neuroimaging data from the ADNI can provide more information regarding diagnosis of MCI than the standardized memory tests obtained as part of the study.
Moreover, we are interested in localizing areas of the cortex that differ between healthy controls (HC) and individuals with MCI.
Three-dimensional T1-weighted structural images for 229 healthy controls and 399 subjects with MCI were obtained as part of the ADNI.
This sample consists of subjects who had images and a composite memory score available at baseline.

We perform the analysis in two ways:
First, we proceed with standard analysis methods currently available for neuroimaging data in open access software \citep{fischl_freesurfer_2012}.
Second, we use the PST statistic and the high-dimensional inference procedure described above.

Cortical thickness was estimated using Freesurfer \citep{fischl_measuring_2000,dale_cortical_1999}. 
Subjects' thickness data were registered to a standard template for analysis and smoothed at 10mm FWHM to reduce noise related to preprocessing and registration.
The template contains 18,715 vertex locations where cortical thickness is measured for each subject.
Our goal is to identify whether the 18,715 cortical thickness measurements provide any additional information regarding the diagnosis of the individuals.
For all analyses we include age, sex, and the composite memory score as covariates \citep{crane_development_2012}.

\subsection{Standard Neuroimaging Analysis Procedure: Average and Vertex-wise Testing}
Because neuroimaging studies typically collect many types of images with many covariates and possible outcomes, it is common to obtain a summary measure of a high-dimensional variable, and then proceed with further analysis if the summary measure appears to be associated with an endpoint of interest.
In this analysis we first take the average of all the cortical thickness measurements across the cortical surface for each subject and perform a regression with diagnosis as the outcome using logistic regression.
Specifically let $C_i$ denote the average cortical thickness measurement for subject $i$, and $X_i$ denote a vector with an intercept term, age, an indicator for sex, and the composite memory score for subject $i$. Then we fit the model
\[
\text{logit}\{\P(Y_i=1 \mid C_i, X_i)\} = X_i^T \alpha + C_i \beta_C.
\]
If there is a significant relationship with the average cortical thickness measurements, i.e. if we reject $H_0: \beta_C = 0$, then we will proceed by performing mass-univariate vertex-wise analyses by running a separate model at each point on the cortical surface.

The analysis using the average cortical thickness variable suggests a highly significant association of cortical thickness with diagnosis, indicating that subjects with thinner cortices are more likely to have MCI (Table \ref{logisticResults}).
Based on these results we choose to investigate the relationship at each vertex to localize where in the cortex the association occurs.

\begin{table}[ht]
\centering
\begin{tabular}{rrr}
  \hline
 & Estimate (SE) & p-value \\ 
  \hline
Age & -0.08 (0.02) & $<0.001$ \\ 
  Sex (Male) & -0.37 (0.26) & 0.15 \\ 
  Memory score & -3.22 (0.27) & $<0.001$ \\ 
  Average cortical thickness & -4.23 (1.02) & $<0.001$ \\ 
   \hline
\end{tabular}
\caption{Results for the logistic regression of diagnosis onto covariates and whole-brain average cortical thickness. Results for average cortical thickness indicate a highly significant association between cortical thickness and diagnosis. SE denotes standard error.} 
\label{logisticResults}
\end{table}

\begin{figure}
  \includegraphics[width=0.98\linewidth]{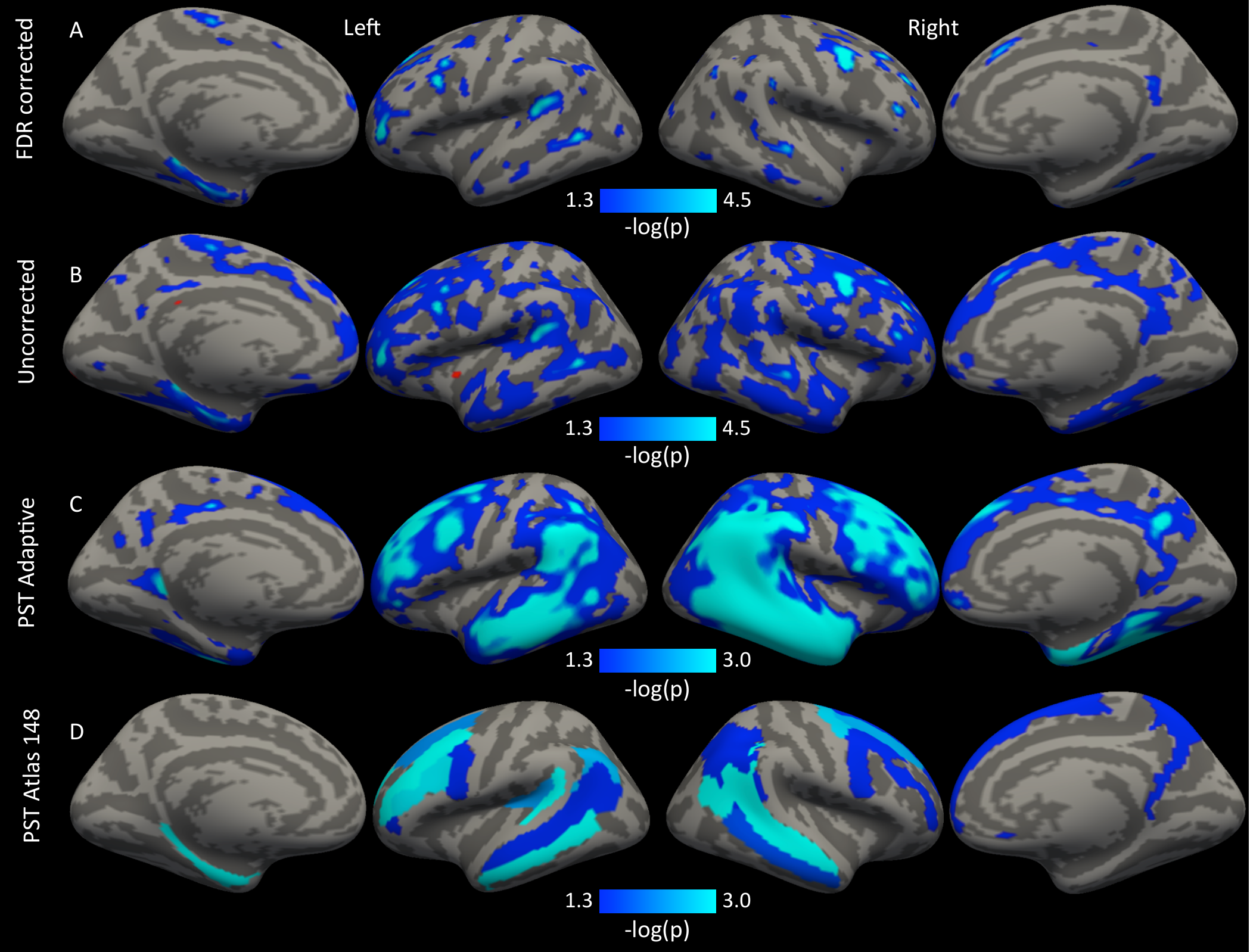}
  \caption{A comparison of inference procedures of the association between the imaging data and diagnosis. (A) Benjamini-Hochberg corrected vertexwise results, (B) Uncorrected vertexwise results, and (C \& D) results based on PST high-dimensional inference that control the FWER of the projected scores. (C) The dimension ($r=7$) of the PCA basis for $\mathbb{L}$ was selected using the adaptive procedure. (D) The 148 dimensional basis constructed from the Destrieux anatomical atlas. Blue values show significant ($\alpha=0.05$) negative association with diagnosis indicating that thinner cortex in these regions is associated with an MCI diagnosis.}
  \label{fig:maxsum20}
\end{figure}

For the vertex-wise analyses, we use the software package Freesurfer to perform Benjamini-Hochberg (BH) correction separately across each hemisphere (Figure \ref{fig:maxsum20} A).
The spatial extent of the FDR-corrected results is more limited than what we might expect given the very strong association between diagnosis and average cortical thickness.
Uncorrected exploratory analyses were conducted to further identify regions related to the whole-brain results (Figure \ref{fig:maxsum20} B).
The most significant results occur in left and right frontal lobes.
These analyses suggest that thinning in larger portions of the frontal and temporal lobes is associated with increased risk of MCI; however, these results are not found using a method that guarantees control of the FWER or FDR.

\subsection{PST and High-Dimensional Inference Procedures}
To use the PST procedure we perform the following steps:
\begin{enumerate}
\item Select a subspace $\mathbb{L}$.
\item Perform the PST for the association between the image and diagnosis.
\item If the test in step 2 rejects, then perform \emph{post hoc} inference as in Section \ref{sec:posthoc}.
\end{enumerate}

We select a basis for $\mathbb{L}$ in the two ways described in Section \ref{sec:specialcases}.
For this analysis we use the aPCA procedure described in Section \ref{sec:adaptPCA} to choose the best PCA basis by testing the first 5 components together and sequentially testing components $6,\ldots,n-m$.
We also present results for the PCA basis fixed at several other dimensions ($r=10,20,50$) to demonstrate how the basis affects the results of the analysis.
In addition we consider a basis constructed from the  $r=148$ regions (74 per hemisphere) of the anatomical atlas of \citet{destrieux_automatic_2010}.
If we were unwilling to condition on the covariance structure of the scores or the anatomical atlas, a basis could be constructed that approximates a predetermined covariance structure (e.g. a spatial AR(1)), or a covariance structure estimated from an independent sample to be used to construct the PCA basis.
In addition to the PST we perform the sequence kernel association test (SKAT) \citep{wu_rare-variant_2011}, the sum of powered scores (SPU) test using the infinity norm, which corresponds to testing the max across the scores \citep{pan_powerful_2014}, and the adaptive sum of powered scores test (aSPU), which has competitive power to many other score tests \citep{pan_powerful_2014}.
The SKAT is known to be more powerful if there is a distributed signal, and the SPU infinity norm will be most powerful for a sparse signal.
The aSPU test combines multiple tests based on the norms $\lVert \S \rVert_\gamma^\gamma$ for $\gamma$ varying over a finite subset of $\mathbb{N}$ by choosing one with the smallest p-value.
Permutation testing is used to assess the significance of these statistics.

The aPCA basis selected $r = 7$ by testing for $r=5$ and then sequentially testing the next two PCs.
With this basis we reject the null hypothesis using the PST (Table \ref{maxsumstatistics}), indicating that there is an association between the image and diagnosis conditioning on the effects of age, sex, and composite memory score.
The test rejects at the $\alpha=0.05$ threshold irrespective of which basis is used.
The SKAT, SPU, and aSPU tests also reject the null.

\begin{table}[ht]
\centering
\begin{tabular}{llrr}
  \hline
 & Test Statistic & p-value & Rejection Threshold \\ 
  \hline
Adaptive (r=7) & 38 & $<$0.001 & 3.2 \\ 
  r=10 & 41 & $<$0.001 & 3.4 \\ 
  r=20 & 50 & $<$0.001 & 3.7 \\ 
  r=50 & 91 & $<$0.001 & 4 \\ 
  Anatomical basis & 179 & 0.04 & 3.4 \\ 
  SKAT & $5\times 10^6$ & $<$0.001 &  \\ 
  SPU Inf & 47 & $<$0.001 &  \\ 
  aSPU &  & 0.001 &  \\ 
   \hline
\end{tabular}
\caption{The $\chi^2$ PST statistic and associated p-values for various basis dimensions; Adaptive, 10, 20, and 50. ``Anatomical" is a basis constructed from an anatomical atlas of dimension 148. The last column denotes the 5\% familywise error rejection thresholds for the projected scores, i.e. the probability any projected score is above those values under the null is 5\%. The thresholds are obtained using 10,000 simulations.} 
\label{maxsumstatistics}
\end{table}

Given the results of the PST we are then interested in investigating how the scores contribute to the significant test statistic.
To investigate the contributions of the scores to the PST statistic we perform \emph{post hoc} inference on the projected scores.
We use 10,000 simulations to obtain rejection regions for each of the basis dimensions.
The simulations ran for all bases in less than 2 minutes.


Results suggest that thinner cortex in bilateral temporal and frontal lobes and right precuneus is associated with an increased risk of MCI (Figure \ref{fig:maxsum20} C \& D).
Results are given as $-\log_{10}(p)$ where $p$ is obtained using the simulated distribution \eqref{eq:pvalues}.
These locations are known to be thinner in AD versus HC as well as in AD versus MCI \citep{singh_spatial_2006} and the results here demonstrate that there are significant differences between MCI and HC in the same region.
The results indicate that the degree of frontal and temporal lobe thinning is correlated with diagnostic severity, and suggest that measurements of cortical thickness may provide useful information over neurospsychological scales in identifying people at risk for AD.
Differences in these regions between MCI and HC were previously shown by \citet{wang_alterations_2009}; however the authors did not control for multiple comparisons or adjust for covariates.

\begin{figure}[h!]
  \includegraphics[width=\linewidth]{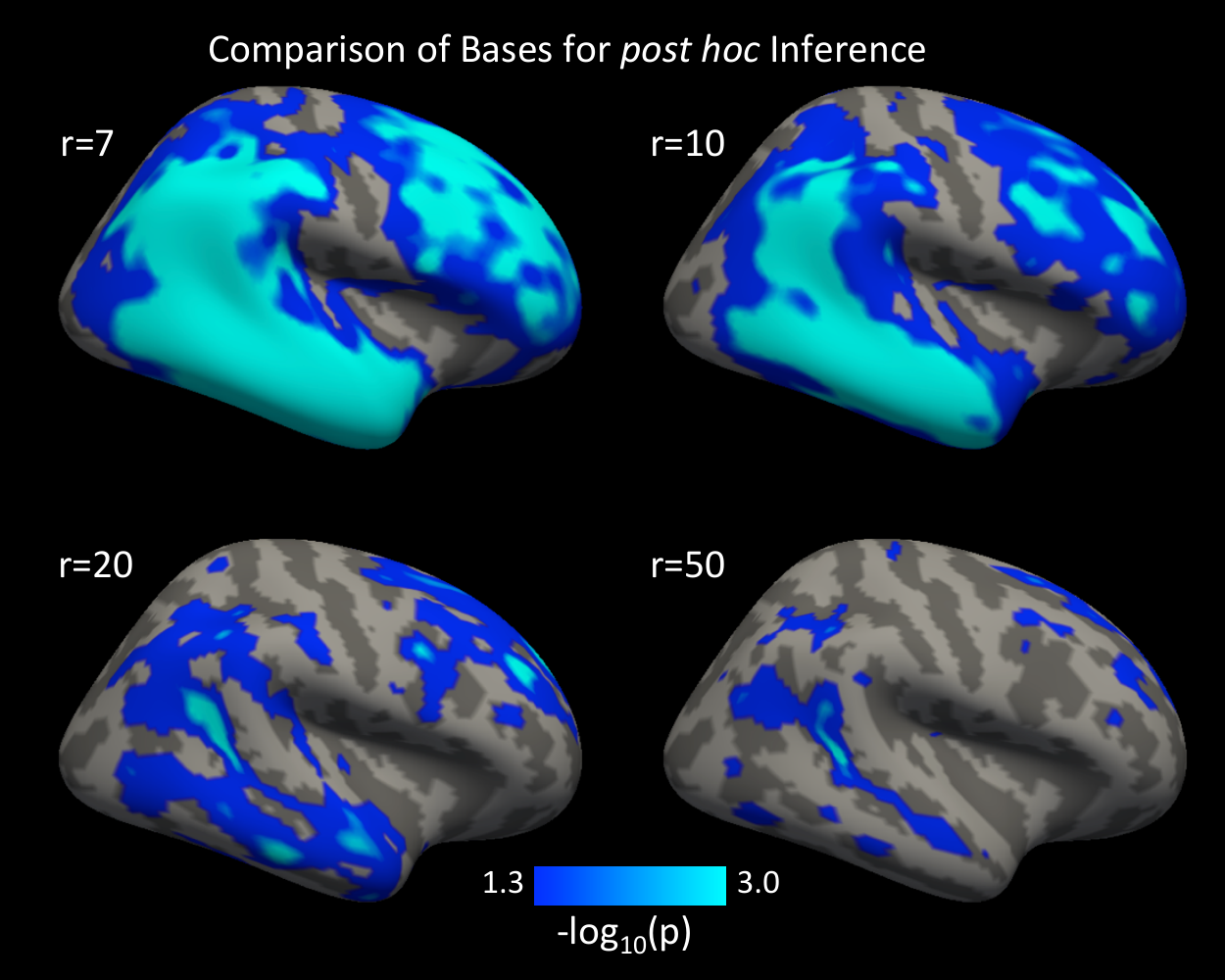}
  \caption{PST inference for PCA bases of various rank; Adaptive (7), 10, 20, 50. Increasing the dimensionality of the basis increases spatial specificity, but comes with the cost of more conservative inference (see e.g. Table \ref{maxsumstatistics}).}
\label{fig:maxsumrank} 
\end{figure}

To reiterate, the blue areas in Figure \ref{fig:maxsum20} C \& D are based on low-rank inference and control the FWER of the projected scores.
The procedure has improved power over standard correction methods seen in Figure \ref{fig:maxsum20} A \& B by performing inference in a lower dimensional space.
The p-values obtained in Figures \ref{fig:maxsum20} and \ref{fig:maxsumrank} use \eqref{eq:pvalues} and indicate the probability of observing a projected score statistic as extreme under the global null $H_0: \beta = 0$.
Though interpretation is restricted to the projected scores, the results align strongly with previous reports \citep{singh_spatial_2006,wang_alterations_2009}.

To demonstrate the impact of the choice of $r$, we performed \emph{post hoc} inference on the scores for 4 different PCA bases (Figure \ref{fig:maxsumrank}).
It is clear from Figure \ref{fig:maxsumrank} that increasing the dimension of the basis increases the spatial specificity of the results. However, the larger bases also come with the cost of reduced power due to the larger degrees of freedom of the basis.
This is also illustrated in Table \ref{maxsumstatistics}, where the larger bases have a higher rejection threshold.

\section{Neuroimaging-based Simulation Study}
\label{sec:Simulations}

As a simulation study, we perform analyses using data generated for the right hemisphere of the cortical thickness data from the ADNI dataset measured at $p=$9,361 locations called vertices.
We simulate an artificial outcome of interest that is categorical, as in the ADNI analyses presented above.
We select two anatomical regions (superior temporal sulcus and superior frontal sulcus) of 669 vertices total to have a negative association with the outcome and one region (anterior part of the cingulate gyrus and sulcus) of 191 vertices to have a positive association.
The first two of these regions were selected because of their association in the ADNI data set.
The third region was selected to compare the performance of the tests when there are different locations with positive and negative associations with the outcome.
To create a mean and covariance structure similar to real data within the regions of association, we create the mean vectors and covariance matrices for the simulations from the full sample of subjects used in the ADNI Freesurfer analysis above, yielding two full rank covariance matrices, $\Sigma_-$ and $\Sigma_+$ and mean vectors $\mu_+$ and $\mu_-$.

For each simulation, we select a random subset of size $n$ without replacement from the subset of control subjects used in the ADNI neuroimaging analysis.
Data within the negatively and positively associated regions are generated as independent multivariate normal distributions for each subject, with covariance structures $G_{i,-} \sim N(\mu_-, \Sigma_-)$ and $G_{i,+} \sim N(\mu_+, \Sigma_+)$, respectively.
We centered the imaging data prior to analysis.

In each simulation the outcome is generated under a logistic model
\begin{equation}\label{eq:simmodel}
\text{logit}(\mathbb{E} Y_i) = \alpha_0 -\beta \mathbf{1}^T G_{i,-} +  2\beta \mathbf{1}^T G_{i,+},
\end{equation}
where $\alpha_0$ is set to the log ratio of MCI to controls in the neuroimaging analyses section.
$\mathbf{1}$, is a vector of ones, and $\beta$ is an unknown parameter that we vary from 0 to 0.005.
We multiply the values in the positive region by 2 to increase signal because it is a spatially smaller cluster than the two negative regions.
In addition to simulations where the coefficients are constant across each region, in the Supplement 
we perform simulations generating the parameters from a uniform distribution.

We construct the subspace $\mathbb{L}$ in three ways. The first is to use the adaptive procedure (Section \ref{sec:adaptPCA}) in each sample conditioning on the estimate $\hat\alpha_0$. The second basis type is constructed in each sample from the first $r = 10, 20, 50$ principal components from a PCA of $G(I-H)$, where $H$ is the projection onto the intercept. The third basis is constructed from regions of anatomical atlas of \citet{destrieux_automatic_2010}, by randomly grouping the 74 regions into $r$ groups and using normed indicator vectors for each group as the basis.

We assess power for indices with nonzero mean and type 1 error for indices with zero mean.
If we denote the set of indices with a nonzero association with the outcome by $J$, then the expectation of the score $\mu_j$ is nonzero only for $j \in J$, where $\mu$ is as defined in \eqref{eq:scoremean}.
So, for indices with $j \notin J$ we report type 1 error and for indices with $j \in J$ we report power.

Similarly, the mean of the standardized projected scores, $\Delta P \mu$, determines type 1 error and power for the projected scores $\Delta P \S$.
The FWER and FDR of the projected scores are reported for the basis constructed from the anatomical atlas and the PCA bases.
In general, no element of the standardized projected mean is exactly zero, so type 1 error is assessed by thresholding the standardized projected parameter vector at the 0.2 quantile and reporting the rejection rate for vertices with projected parameter values below that threshold.


We perform 1000 simulations for sample sizes of $n=100, 200$ and compare the PST for the adaptive procedure and fixed bases with dimensions of $r=10, 20, 50$.
In addition, we compare the PST to the sequence kernel association test (SKAT) \citep{wu_rare-variant_2011}, the sum of powered scores (SPU) test using the infinity norm, and the adaptive sum of powered scores test (aSPU) \citep{pan_powerful_2014}.
We assess pointwise power and type 1 error of the PST inference with uncorrected, Bonferroni-corrected, and BH-corrected results.
We also compare FWER and FDR between methods. For these comparisons we assess the type 1 error for the unprojected scores using inference designed for the projected scores.

\begin{figure}
\center{
  \includegraphics[width=1\linewidth]{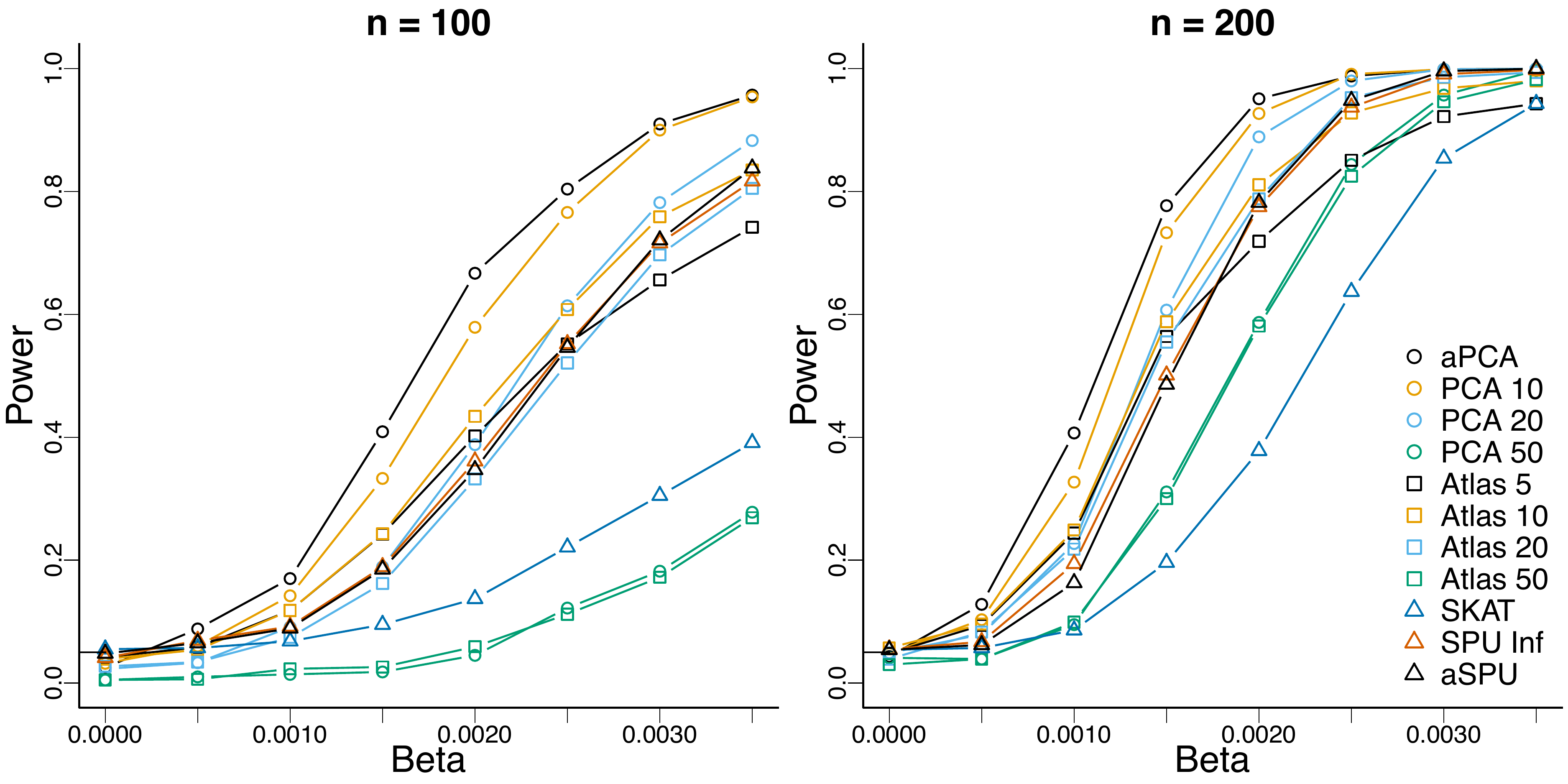}}
  \caption{Power results for the PST with various bases compared to aSPU and the SKAT. ``aPCA" is the adaptive basis and ``PCA $r$" indicates the basis formed from the first $r$ components of the PCA of the design matrix $G$. ``Atlas $r$" indicates the bases formed from the anatomical atlas with $r$ regions covering the cortex.}
\label{fig:maxsumpower} 
\end{figure}

The PST with fixed bases demonstrates superior power to the other tests (Figure \ref{fig:maxsumpower}), due to its ability to remove the influence of unassociated scores from the test by maximizing over the basis and by leveraging the spatial information in the data.
If these features of the data were not informative then the PST would not perform well. 
The aPCA is has better power than the other PCA bases because a low rank basis is all that is required to capture the signal in the data.
aSPU is adaptive to the sparsity of the signal, so it performs better than the SKAT, but does not use the information in the covariance of the scores to leverage power.

\begin{table}[ht]
\centering
\begin{tabular}{lrr}
  \hline
 & FWER & FDR \\ 
  \hline
aPCA & 0.06 & $<0.01$ \\ 
  PCA 10 & 0.04 &$<0.01$ \\ 
  PCA 20 & 0.03 & $<0.01$ \\ 
  PCA 50 & 0.02 & $<0.01$ \\ 
  Anatomical 5 & 0.03 & 0.02 \\ 
  Anatomical 10 & 0.04 & 0.02 \\ 
  Anatomical 20 & 0.04 & 0.02 \\ 
  Anatomical 50 & 0.06 & 0.02 \\ 
   \hline
\end{tabular}
\caption{Error rates for the projected scores for the adaptive PCA bases and anatomical bases for $n=200$ and $\beta=0.002$.}
 \label{tab:projerror}
\end{table}
As expected, the \emph{post hoc} inference procedure controls the FWER of the projected scores for all basis dimensions (Table \ref{tab:projerror}).
In general, the \emph{post hoc} inference procedure does not control the FWER or FDR of the unprojected scores (Table \ref{tab:unprojerror}) as the inference is intended for the projected scores.
However, for larger PCA bases our procedure does control the FDR (bold rows in Table \ref{tab:unprojerror}).
This is likely because the projection captures most of the variation in $\mu$, so that the projection $\Delta P\mu$ is close to $\mu$.
Future investigation of whether inference for the projected scores will control any error rate for unprojected score vector is warranted.

\begin{table}[ht]
\centering
\begin{tabular}{lrrr}
  \hline
 & HR & FDR & FWER \\ 
    \hline
  aPCA & 0.69 & 0.13 & 0.70 \\ 
  \bf{PCA 10} & \bf{0.55} & \bf{0.05} & \bf{0.45} \\ 
  \bf{PCA 20} & \bf{0.34} & \bf{0.02} & \bf{0.23} \\ 
  \bf{PCA 50} & \bf{0.12} & \bf{0.01} & \bf{0.10} \\ 
  Anatomical 5 & 0.16 & 0.79 & 0.27 \\ 
  Anatomical 10 & 0.34 & 0.64 & 0.57 \\ 
  Anatomical 20 & 0.53 & 0.45 & 0.84 \\ 
  Anatomical 50 & 0.63 & 0.17 & 0.70 \\ 
  Uncorrected & 0.61 & 0.44 & 1.00 \\ 
  \bf{Holm} & $\mathbf{<0.01}$ & $\mathbf{<0.01}$ & $\mathbf{<0.01}$ \\ 
  \bf{BH} & \bf{0.10} & \bf{0.05} & \bf{0.22} \\ 
\end{tabular}
\caption{Hit (HR), false discovery (FDR), and family-wise error (FWER) rates for the unprojected scores. Higher dimensions of the PCA basis control the FDR while maintaining a higher hit rate for the unprojected scores. Note, however, that in practice it is not possible to know the dimension of the basis required to control the FDR of the unprojected scores. Bold rows control the FDR at $q=0.05$. BH=Benjamini-Hochberg.}
\label{tab:unprojerror}
\end{table}

\begin{figure}
\center{
  \includegraphics[width=1\linewidth]{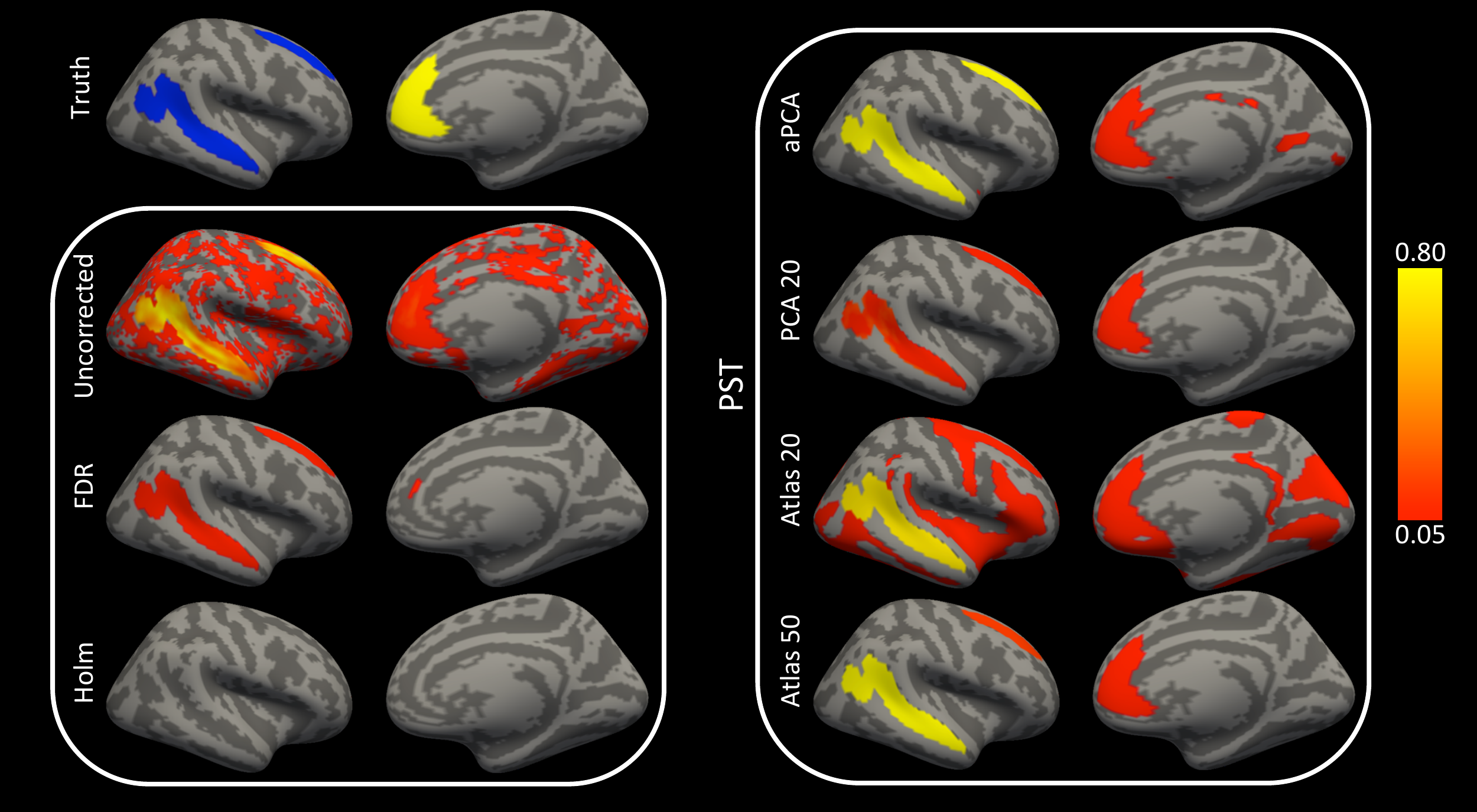}}
  \caption{Vertexwise power or type 1 error is measured as the proportion of simulated samples where each of the testing procedures rejects the null at a given location. Results are shown for $\beta = 0.002$. ``Truth" indicates locations where signal was simulated according to the model \eqref{eq:simmodel}.}
\label{fig:maxsumvertexwisepower} 
\end{figure}

The vertexwise error rate describes how effective a procedure is at controlling the error rate for the unprojected scores at each location.
The vertexwise error rate of the PST inference procedure for the unprojected scores using the PCA 10 basis is low while maintaining better vertexwise power than BH (Figure \ref{fig:maxsumvertexwisepower}; PCA 10). This is because in any given sample there may be a high false positive rate, but the errors across samples do not appear in the same locations. The BH and Bonferroni corrections both work well at controlling the vertexwise type 1 error rate but have lower power compared to the PCA-based PST procedure (Figure \ref{fig:maxsumvertexwisepower}).
The bases constructed from the anatomical atlas tend to have large regions of vertexwise type 1 error for the unprojected scores. At the largest basis dimension the atlas allows for enough specificity to reduce the vertexwise error.
All methods have lower power to detect the positive cluster than the two negative clusters.
This is possibly due to the characteristics of the covariance structure in the positive cluster which overlaps gyral and sulcal regions.

\section{Discussion}
We have proposed the PST, which maximizes the weights for the Sum statistic in a subspace of the parameter space.
The procedure offers a novel \emph{post hoc} inference on the projected scores by performing inference in the subspace where the test statistic was estimated.
Because the posthoc inference is based on the same model and degrees of freedom as the PST statistic, the interpretation of high-dimensional results agree closely with the results from the PST.

Instead of choosing a specific value for the weight vector, $\zeta$ as in the sum test, our methodology allows the investigator to select a space to consider for $\zeta$.
The ability to choose a space makes the procedure very flexible.
For example, in imaging the basis for the space can be chosen based on anatomical or functional labels, or from data acquired in another imaging modality.
Particular hypotheses can be targeted by selecting a basis that includes indicators of certain regions or weights particular locations to target specific spatial patterns.
If orthogonal indicator vectors are used as the basis, then the approach can be seen as testing averages of subregions of the data as in Section \ref{sec:DataAnalysis}.
In this case, the PST procedure can be seen as a maxT multiple testing procedure that accounts for the correlation structure of the tests.

There are several limitations of the proposed procedure.
First, the success of the procedure depends critically on the projection chosen.
If a projection is chosen that is orthogonal to the mean vector then the PST will fail to capture any signal in the data.
This is a limitation of any dimension-reducing procedure.
Further research could investigate whether maximization of the score test with regularization can yield a test statistic whose distribution is tractable.
Regularization may remove the subjectivity of selecting a basis and make the procedure more robust.
Second, while the dimension reduction procedure preserves power and the results align closely with those from previous research, the inference does not guarantee control of the FWER or FDR of the original score vector.
Future research will investigate how inference of the original score vector can be made by thresholding the projected score vector.
This is similar in concept to the dependence-adjusted procedure discussed by \citet{fan_estimating_2012} for controlling the FDP and may offer increased power by leveraging the covariance of the test statistics.
These limitations notwithstanding, our procedure generalizes Rao's score test to the high- and infinite- dimensional settings and introduces a new inference approach based on projecting the test statistics to a lower-dimensional space where inference can be made on fewer degrees of freedom.


\appendix

\makeatletter   
 \renewcommand{\@seccntformat}[1]{{\csname the#1\endcsname}.\hspace*{1em}}
 \makeatother

\renewcommand{\theequation}{A\arabic{equation}}
\setcounter{equation}{0}
\renewcommand{\thefigure}{A\arabic{figure}}
\setcounter{figure}{0}
\renewcommand{\thetable}{A\arabic{table}}
\setcounter{table}{0}
 
\section{APPENDIX A}

\subsection{Theoretical framework}
\label{sec:framework}
We assume the observed data are finite-dimensional representations that are generated from an underlying stochastic process.
To be more specific, define the Hilbert space $\mathbb{Y} =\R^k\times \mathcal{B}(\mathbb{V}) $, where $\mathbb{V}$ is a nonempty compact subset of $\R^3$ and $\mathcal{B}(\mathbb{V})$ is the space of square integrable functions from $\mathbb{V}$ to $\R$.
$\mathbb{V}$ represents the space on which data can be observed; in neuroimaging this space is the volume of the brain.
Let $(\mathbb{Y}, \mathcal{Y}, \P)$ be a probability space where $\mathcal{Y}$ is the Borel $\sigma$-algebra on $\mathbb{Y}$.
Let $\Y_i = (\Y_i^{(1)}, \Y_i^{(2)})$, for $i=1,\ldots, n$, be iid with $\Y_i^{(1)}$ taking values in $\R^k$, $\Y_i^{(2)}$  a stochastic process taking values in $\mathcal{B}(\mathbb{V})$, and $\P(\Y_i \in A) = \P(A)$ for all $A \in \mathcal{Y}$.
Observations of $\Y_i$ are bounded functions together with a vector of $k$ variables that are nonimaging covariates and the outcome variable.
We denote the collection of data by $\Y = (\Y_1, \ldots, \Y_n)$.
Although $\Y$ represents the underlying data, in practice $\Y_i^{(2)}$ are unobservable and we only observe discretized data at a finite number of locations that are voxels in the image.

We define a parameter space $\Theta = \Alpha \times \Beta$ that includes a finite-dimensional parameter $\alpha \in \Alpha \subset \R^{m}$ and an infinite-dimensional parameter on $\beta \in \Beta$.
Together these parameters describe the joint distribution of the imaging and nonimaging data.
Throughout, we assume $\Beta = \mathcal{B}(\mathbb{V})$, so that the infinite-dimensional parameter and infinite-dimensional data are defined on the same space, but this assumption is not required.
The distribution of the observed data will be defined by a $p$-dimensional discretization of the infinite-dimensional parameter.
We will prove that under a few assumptions, as $n,p \to \infty$ the test statistic for the discretized data approaches the statistic for the infinite-dimensional parameter.

To relate the unobserved data $\Y$ to the parameters, we further assume that the measure $\P$ is in a set of probability models, $\{ F_\theta : \theta \in \Theta \}$, indexed by the parameter 
\begin{equation}\label{eq:parameterdef}
\theta = (\alpha, \beta).
\end{equation}
That is, there exists a regular point $\theta_0 \in \Theta$ such that for all sets $A \in \mathcal{Y}$
\[
\P(\Y_i \in A) = F_{\theta_0}(A).
\]
We define the density function $f_\theta$ as the Radon-Nikodym derivative of $F_\theta$ with respect to the Lebesgue measure $\mu$,
\[
\P\left( \Y_i \in A \right )
= \int_A \frac{dF_\theta}{d\mu} d\mu
= \int_A f_\theta d\mu.
\]
Let $\ell(\theta; \Y) = n^{-1}\sum_{i=1}^n\log f_\theta(\Y_i)$ be the log-likelihood function for $\theta$.

In order to define the finite-dimensional data and parameter space we must partition $\mathbb{V}$ into finitely many sets and define the observable random variables as a realization from the partitioned space.
For any integer $p$, the space $\mathbb{V}$ can be partitioned into $p$ nonempty sets. 
Denote the partition $\mathcal{V}_p = \{ \mathbb{V}_{1p}, \ldots, \mathbb{V}_{pp}\}$. Note that, by the definition of a partition, $\bigcup_{j=1}^p \mathbb{V}_{jp} = \mathbb{V}$ and $\mathbb{V}_{jp} \cap \mathbb{V}_{kp} = \varnothing$.
Let $v_j$ be an arbitrary interior point of $\mathbb{V}_{jp}\in \mathcal{V}_p$.
Let the discretized data be $\Y_{ip} = (\Y_{i}^{(1)}, \Y_{i}^{(2)}(v_1), \ldots, \Y_{i}^{(2)}(v_p)) \in \R^{k+p}$ whose distribution is determined by the finite parameter $\theta_p = (\alpha, \beta(v_1), \ldots, \beta(v_p)) \in \R^{m+p}$.
In order to define a finite-dimensional likelihood from the likelihood for the infinite-dimensional parameters we define the function $\beta_p \in \mathcal{B}(\mathbb{V})$ by
\[
\beta_p(v) = \beta(v_j) \text{ for } v \in \mathbb{V}_{jp}
\]
and the stochastic processes
\[
\Y^{(2)}_{ip}(v) = \Y^{(2)}_i(v_j) \text{ for } v \in \mathbb{V}_{jp}.
\]
As $\mathcal{V}_p$ is a partition, each $v$ is in only one $\mathbb{V}_{jp}$.
This allows us to define the log-likelihood from the function
\begin{equation}
\label{eq:finitelikelihood}
\ell(\theta_p; \Y_{p}) = n^{-1} \sum_{i=1}^n \log f_{\theta_p}(\Y_{ip}),
\end{equation}
where $f_{\theta_p}(y_{p}) = f((y^{(1)}, y_{p}^{(2)} ); (\alpha, \beta_p) )$.

Finally, assuming $\ell$ is Fr\'echet differentiable with respect to $\beta$, we can define the scores
\begin{align}
\U_n = \U_n(v) &= \frac{\partial \ell}{\partial \beta}\{(\alpha, \beta(v)); \Y(v)\}\label{eq:infinitederivfunc}\\
\U_{np} = \U_{np}(v) &= \frac{\partial \ell}{\partial \beta_p}\{(\alpha, \beta_{p}(v)); \Y_p(v)\}\label{eq:frechetfinitederivfunc} .
\end{align}
Let 
\begin{align*}
\mathrm{S}_{n} & = \mathrm{U}_{n}(\cdot ; (\hat\alpha, \beta_0) ) \in \mathcal{B}(\mathbb{V}) \\
\mathrm{S}_{np} & = \mathrm{U}_{np}(\hat\alpha, \beta_{p0})  \in \R^p
\end{align*}
where $\beta_0$ denotes the value of the parameter under the null  $H_0: \beta = \beta_0$ and $\hat \alpha$ is the maximum likelihood estimator for $\alpha$ under the null.


\subsection{Conditions for Theorem \ref{thm:MaxSumGeneral}}
\label{conditions}
 The conclusion of Theorem \ref{thm:MaxSumGeneral} requires the asymptotic normality of the scores, which holds under the following conditions:
 \begin{enumerate}
 \item The ability to interchange integration and expectation of the likelihood so that 
 \[
 \E_{\theta_{p0}} \S_{np}
 = \frac{\partial}{\partial \beta_p} \E_{\theta_{p0}} \ell\{(\alpha, \beta_{p}(v)); \Y_p(v)\}
 = 0.
 \]
 \item The score  \eqref{eq:effinfo}  variance is finite.
 
 Asymptotic normality of the scores then follows by the multivariate central limit theorem \citep{van_der_vaart_asymptotic_2000}.
 
 \end{enumerate}

 \subsection{Proof of Theorem \ref{thm:MaxSumGeneral}}
 \label{prooftheorem1}

 Let $\phi = Q^T \zeta$. Then the PST statistic is
\begin{align}
\mathrm{R}^{\mathbb{L}}
& = \max_{\zeta \in \R^p\setminus \{0\}}n \frac{(\zeta^T P \S_{np})^2}{\zeta^T P \hat\Omega P \zeta } \notag \\
& = \max_{\phi \in \R^r\setminus \{0\}}n \frac{(\phi^T Q^T \S_{np})^2}{\phi^T \hat{\mathrm{V}} \phi } \notag \\
& = n \mathrm{S}_{np}^T Q \hat{\mathrm{V}}^{-1} Q^T \mathrm{S}_{np} \notag,
\end{align}
where the last line follows from a standard maximization lemma \citep[p.~80]{johnson_applied_2007}.

Equation \eqref{eq:rotatedscores} holds by the multivariate central limit theorem and the variance estimate of $n^{1/2}Q^T \S_{np}$ is
\begin{align*}
\hat{\mathrm{V}}(\theta_0)
& = Q^T\hat \Omega(\theta_0) Q \\
& = (Q^T\hat \Omega_\beta Q - Q^T\hat\Omega_{\beta\alpha}\hat\Omega^{-1}_{\alpha} \hat\Omega_{\alpha\beta} Q),
\end{align*}

which converges to  $V(\theta_0)$ by the continuous mapping theorem because $\hat\Omega_F \to_P \Omega_F$.
Thus $n \mathrm{S}_{np}^T Q \hat{\mathrm{V}}^{-1} Q^T \mathrm{S}_{np} \to_L \chi^2_r$ by the continuous mapping theorem.

\begin{remark}
\normalfont The conclusion of Theorem \ref{thm:MaxSumGeneral} implies that expression \eqref{thm1maxsum} does not depend on the choice of $Q$. This fact can also be shown directly, as follows. Consider another matrix $Q_*$ with orthonormal columns such that $P=Q_*Q_*^T$, and accordingly define $\hat{V}_*=Q_*^T\hat{\Omega}(\theta_0)Q_*$. Then $Q_*=QM$ where $M= Q^TQ_*$. Since $P=QMQ_*^T$ is of rank $r$, $M$ is of rank $r$ and hence invertible, so
\[Q_*\hat{V}_*^{-1}Q_*^T=QM(M^T\hat{V}M)^{-1}M^TQ^T=Q\hat{V}^{-1}Q^T,\]
and thus formula \eqref{thm1maxsum} for the PST statistic is unchanged by substituting $Q_*,\hat{V}_*$ for  $Q,\hat{V}$.
\end{remark}

\subsection{Details for Section \ref{sec:asympp}}
\label{sec:extras_asympp}

For  $y = (y^{(1)}, \ldots, y^{(k)}) \in \mathbb{D}_1\times \ldots \times \mathbb{D}_k$, where $\mathbb{D}_j$ are Hilbert spaces, we define the norm
\[
\lVert y \rVert = \sup_{j} \lVert y^{(j)} \rVert.
\]
Then, following \citet{van_der_vaart_asymptotic_2000}, define a derivative on the space $\Theta$ from Section \ref{sec:MaxSumExponential} .
\begin{definition}
For $\Theta$ as defined in Section \ref{sec:MaxSumExponential}, a function $f:\Theta \to \R$ is called Fr\'echet differentiable at $\theta$ if there exists a bounded linear operator $A: \Theta \to \R$ such that
\[
\lim_{\lVert h\rVert \to 0} \frac{\lVert f(\theta+h) - f(\theta) - Ah \rVert}{ \lVert h \rVert},
\]
for $h \in \Theta$. 
\end{definition}

The following theorem \citep[Thm 18.14]{van_der_vaart_asymptotic_2000} gives conditions under which $n^{1/2}\mathrm{S}_n \to_L \mathrm{S}$, where $\mathrm{S}$ is a mean zero Gaussian process.
\begin{theorem}
\label{thm:stochasticprocessCLT}
The sequence of elements $\sqrt n \mathrm{S}_n$ converges in law to a mean zero Gaussian process $\mathrm{S}$ if and only if
\renewcommand{\theenumi}{\alph{enumi}}
\begin{enumerate}
\item The sequence $n^{1/2}(\S_{n}(v_1), \ldots, \S_n(v_p))$ converges in distribution in $\R^p$ for every finite set of points $v_1,\ldots, v_p \in \mathbb{V}$.
\item For every $\epsilon, \eta >0$ there exists a partition of $\mathbb{V}$ into finitely many sets $\mathbb{V}_1, \ldots, \mathbb{V}_p$ such that
\[
\limsup_{n\to\infty}\P \left( \sup_j \sup_{ v_1, v_2 \in \mathbb{V}_j} \lvert \S_{n}(v_1) - \S_n(v_2) \rvert \ge \epsilon \right) \le \eta.
\]
\end{enumerate}
\renewcommand{\theenumi}{\arabic{enumi}}
\end{theorem}
Condition (a) of Theorem \ref{thm:stochasticprocessCLT} is satisfied under the assumptions in Appendix \ref{conditions}.
Condition (b) implies that the process $\S$ is continuous in probability.
We require this as an assumption for asymptotics in $p$.

\subsection{Proof of Theorem \ref{thm:maxsumconvergence}}
\label{sec:proofasympp}

If we show
\begin{equation}\label{eq:limits}
\lim_n \lim_p \sqrt n Q_p^T \S_{np} =_L \lim_p \lim_n \sqrt n Q_p^T \S_{np} =_L Q^T \S,
\end{equation}
 where $=_L$ denotes equality in distribution and
\begin{equation*}
\hat{\mathrm{V}}_{np} \to_P V,
\end{equation*}
then the continuous mapping theorem implies \eqref{eq:maxsumconvergence}.
Expanding an arbitrary element of the vector on the left-hand side of \eqref{eq:limits} gives
\begin{align*}
\lim_n \lim_p \sqrt n (Q_p^T \S_{np})_j
& = \lim_n \lim_p \sqrt n \sum_{k=1}^p q_{j}(v_{kp}) \S_n(v_{kp}) \nu(\mathbb{V}_{kp}) \\
& = \lim_n \sqrt n \int_\mathbb{V} q_j(v) \S_n(v)d\nu(v) \\
& =_D \int_{\mathbb{V}} q_j(v) \S(v) d \nu(v) = (Q^T \S)_j.
\end{align*}
The second equality follows because the numerator on the right-hand side is a Riemann integral and our assumptions \eqref{eq:discontinuities} and that $\S_n(\cdot ; \Y)$ is continuous for almost all $\Y$ in Theorem \ref{thm:maxsumconvergence} guarantee that $q_j(v) \S_n(v)$ is integrable.
The final equality follows from the continuous mapping theorem since $\S_n \to_L \S$.
For the limit
\begin{align*}
\lim_p \lim_n \sqrt{n}\sum_{k=1}^p q_{j}(v_{kp}) \S_n(v_{kp}) \nu(\mathbb{V}_k)
& =_L \lim_p \sum_{k=1}^p q_{j}(v_{kp}) \S(v_{kp}) \nu(\mathbb{V}_k) \\
& = \int_{\mathbb{V}} q_j(v) \S(v) d \nu(v) = (Q^T \S)_j,
\end{align*}
where the first line applies the continuous mapping theorem to the finite dimensional vector $\S_{np}$ and the continuity of $\S$ implies the integral exists.
For both directions the limit on $p$ requires \eqref{eq:intconvergence}, so that the volume of all voxels goes to zero.
Theorems \ref{thm:MaxSumGeneral} and \ref{thm:stochasticprocessCLT} are needed to ensure that $\S_{np} \to \S_p$ and $\S_n \to \S$.
The proof for the convergence of $\hat{\mathrm{V}}_{np} \to_P V$ is a similar argument and relies on the assumption that the sample paths of $\partial/\partial \theta \log f(\Y_i; \theta(v))$ are continuous almost everywhere, so that the Riemann integral converges.

 \subsection{Proof of Theorem  \ref{thm:MaxSum}}
 \label{prooftheorem2}

We will ignore the term $\frac{1}{\hat \sigma^2}$ in the maximization as it is constant with respect to $\zeta$.
Define $\mathbb{M} = G \mathbb{L}$ to be the column space of $W$.
From the definition of $\mathrm{R}^{\mathbb{L}}$
\begin{align}
\mathrm{R}^{\mathbb{L}} \propto &
\max_{\zeta \in \mathbb{R} ^p}\frac{(\Y^T G P \zeta)^2}{ \zeta^T P G^T G P \zeta }\notag \\
& = \max_{\zeta \in \mathbb{R}^p }\frac{(\Y^T W G P \zeta)^2}{ \zeta^T P G^T W G P \zeta } \notag \\
& = \max_{\gamma \in \mathbb{M} }\frac{(\Y^T W \gamma)^2}{ \gamma^T W \gamma }, \notag \\
& = \max_{\gamma \in \mathbb{R}^p }\frac{(\Y^T W \gamma)^2}{ \gamma^T W \gamma },\label{eigenMaxSum}
\end{align}
where $\gamma =G P \zeta$.
Because $\gamma$ must be in the subspace $\mathbb{M}$ since it is the column space of $W$ and $GP$.
The solution to the Rayleigh quotient \eqref{eigenMaxSum} is the solution to the largest generalized eigenvalue problem,
\[
W \Y \Y^T W \gamma = \lambda_{max} W \gamma,
\]
where $\lambda_{max} \in \R$, and $\lVert \gamma \rVert = 1$.
By letting $\phi = W\gamma$ the solution is equivalent to the largest eigenvalue problem
\[
W\Y \Y^T W \phi = \lambda_{max} \phi
\]
Then
\[
\text{eigmax}( W\Y \Y^T W) = \text{eigmax}( \Y^T W \Y) = \Y^T W \Y.
\]
Thus, we have
\[
\mathrm{R}^{\mathbb{L}} = (n-m)\frac{\Y^T W \Y}{\Y^T \Y}.
\]
To derive \eqref{eq:normalmaxsumdist}, note
\begin{equation}\label{scheffe}
\frac{\Y^T \Y}{\Y^T W \Y} = 1 + \frac{\Y^T (I-W) \Y}{\Y^T W \Y}.
\end{equation}
The numerator and denominator of the random term on the left hand side are independent since $P(I-P) = 0$.
So \eqref{scheffe} is distributed as $ 1 + \frac{n-m-r}{r} F_{(n-m-r),r}$.
Thus
\begin{align*}
\mathrm{R}^{\mathbb{L}} 
& =_L \frac{(n-m)}{1 + \frac{(n-m-r)}{r}F_{(n-m-r),r}} \\
& = \frac{r(n-m)}{r + (n-m-r)F_{(n-m-r),r}}.
\end{align*}

\subsection{The \emph{post hoc} inference procedure controls the FWER}
\label{sec:subsetpivotality}

The assumption of subset pivotality states that for any set $I\subset \{ 1,\ldots, p \}= H$ the distribution of the maximum of the standardized projected scores \eqref{eq:RR} in set $I$, given indices in $I$ are true nulls is equal to the distribution of the maximum given all hypotheses are true nulls \citep{westfall_multiple_2008}
\[
\max_{j\in I}\left\{  \lvert (\Delta P \S)_j \rvert \mid \text{$I$ are null} \right\}= \max_{j\in I} \left\{ \lvert (\Delta P \S)_j \rvert \mid \text{$H$ are null}\right\}.
\]
By true null we mean that the expectation of the projected score is zero.
This assumption allows us to construct rejection regions assuming all scores are true nulls, but still have strong control of the FWER, i.e. in the case that $\Delta P \mu\ne 0$, where $\mu$ is defined in \eqref{eq:scoremean}.
Subset pivotality is satisfied for normally distributed statistics because the covariance of the statistics is not affected by changing the mean structure.
Thus, the maximum over a subset of true nulls is not affected by the value of the mean for the other statistics.
So long as the covariance estimates for the true nulls are consistent we will maintain asymptotic control of the FWER.
In our \emph{post hoc} inference procedure the variance estimates are consistent because 
\begin{align*}
(\Delta P \hat\Omega P\Delta)_{jj} & \to_P \text{Var}\{(\Delta P \S)_j\} + (\Delta P\mu)_j \\
(\Delta P \hat\Omega P\Delta)_{jk} & \to_P \E\{ (\Delta P \S)_j (\Delta P \S)_k\}.
\end{align*}
Because $(P\mu)_j = 0$ for all true nulls the variance and covariance estimates are consistent for all of true nulls.


 
\section{Supplementary Material}

\subsection{Simulation analyses}
\label{app:simulations}

In addition to the simulations performed in Section \ref{sec:Simulations} we performed similar simulations in which the data were generated from the model
\begin{equation*}
\text{logit}(\mathbb{E} Y_i) = \alpha_0 -\beta \omega_1^T G_{i,-} +  2\beta \omega_2^T G_{i,+},
\end{equation*}
where $\omega_{1j} \sim \text{Unif}(0.5,1.5)$ and  $\omega_{2j} \sim \text{Unif}(1, 3)$.

\begin{figure}
\center{
  \includegraphics[width=1\linewidth]{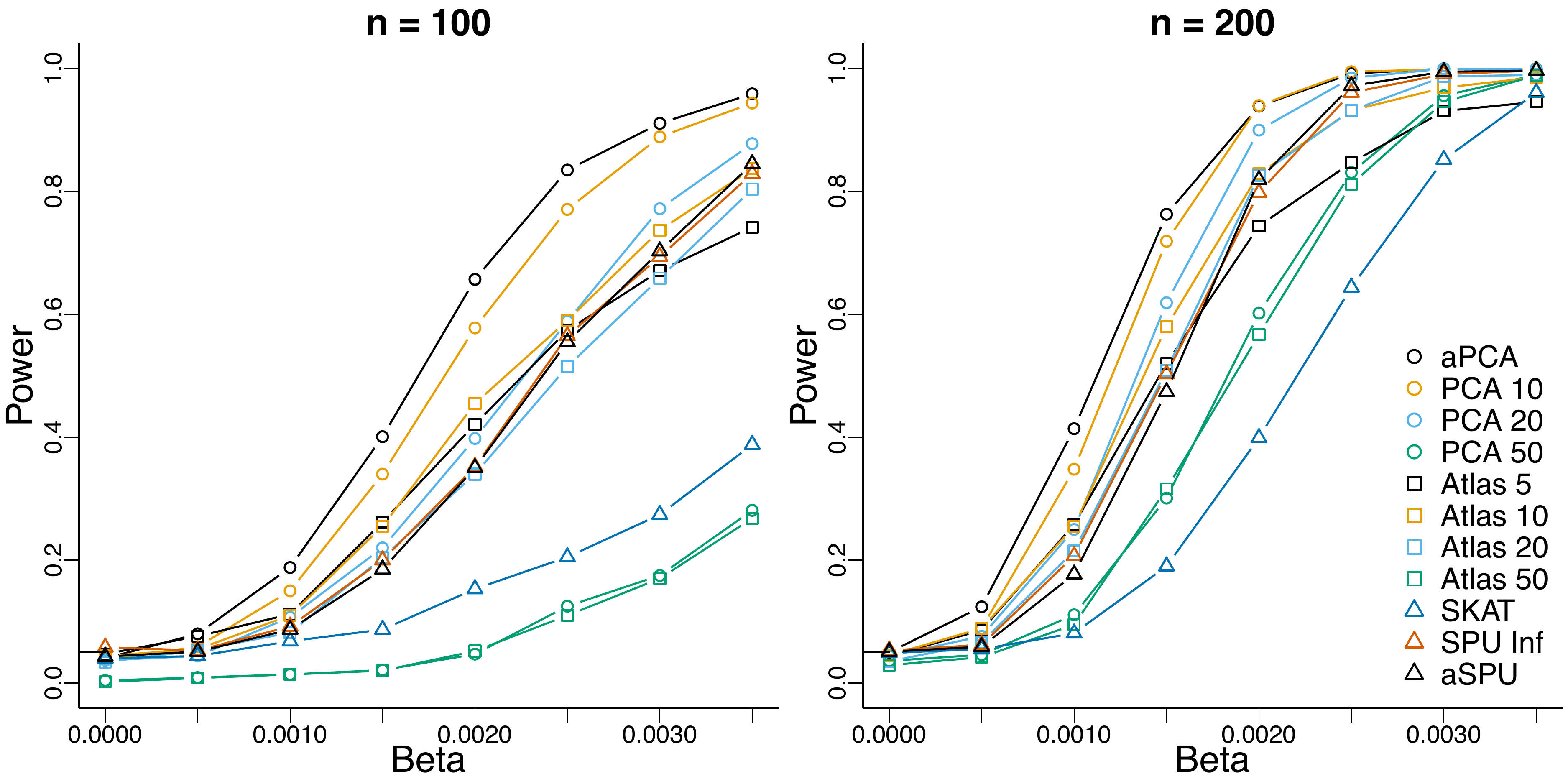}}
  \caption{Supplementary power results for the PST with various bases compared to aSPU and the SKAT. ``aPCA" denotes the adaptive PCA procedure. ``PCA $r$" indicates the basis formed from the first $r$ components of the PCA of the design matrix $G$. ``Atlas $r$" indicates the bases formed from the anatomical atlas with $r$ regions across the cortex. Results are similar to those from simulations performed with constant coefficients.}
\label{suppfig:maxsumpower} 
\end{figure}

Table \ref{supptab:projerror} gives FWER and FDR for the simulation analyses with uniform coefficients.
Figure \ref{suppfig:maxsumpower} gives power results for the additional simulations.
Table \ref{supptab:projerror2} gives FWER and FDR for the simulations where there is no association between image and outcome.

\begin{table}[h!]
\centering
\begin{tabular}{lrr}
  \hline
 & FWER & FDR \\ 
  \hline
aPCA & 0.06 & $<0.01$ \\ 
  PCA 10 & 0.04 & $<0.01$ \\ 
  PCA 20 & 0.02 & $<0.01$ \\ 
  PCA 50 & 0.02 & $<0.01$ \\ 
  Anatomical 5 & 0.04 & 0.02 \\ 
  Anatomical 10 & 0.04 & 0.02 \\ 
  Anatomical 20 & 0.04 & 0.02 \\ 
  Anatomical 50 & 0.05 & 0.02 \\ 
   \hline
\end{tabular}
\caption{Supplementary simulation analysis error rates for the projected scores for the anatomical bases of dimension 5, 10, 20 and 50 for $n=200$ and $\beta=0.002$.}
 \label{supptab:projerror}
\end{table}

\begin{table}[h!]
\centering
\begin{tabular}{rrr}
  \hline
 & FWER & FDR \\ 
  \hline
aPCA & 0.03 & 0.01 \\ 
  PCA10 & 0.03 & 0.01 \\ 
  PCA20 & 0.02 & 0.01 \\ 
  PCA50 & 0.02 & 0.01 \\ 
  Destrieux5 & 0.03 & 0.02 \\ 
  Destrieux10 & 0.04 & 0.03 \\ 
  Destrieux20 & 0.04 & 0.04 \\ 
  Destrieux50 & 0.05 & 0.05 \\ 
   \hline
\end{tabular}
\caption{Supplementary simulation analysis error rates for the projected scores for the PCA and anatomical bases for $n=200$ and $\beta=0$.}
 \label{supptab:projerror2}
\end{table}

\bibliographystyle{plainnat}
\bibliography{MyLibrary}

\end{document}